 \documentclass[11pt]{article}

\usepackage{amssymb}

\newtheorem{theorem}{$~~~~$ Theorem}[section]

\newtheorem{example}[theorem]{$~~~~$ Example}

\newtheorem{corollary}[theorem]{$~~~~$ Corollary}
\newtheorem{lemma}[theorem]{$~~~~$ Lemma}

\newtheorem{definition}[theorem]{$~~~~$Definition}

\newtheorem{dff}[theorem]{$~~~~$ Definition}
\newtheorem{exx}[theorem]{$~~~~$ Example}

\newtheorem{propp}[theorem]{$~~~~$ Proposition}
\newtheorem{remm}[theorem]{$~~~~$Remark}
\newtheorem{ccr}[theorem]{$~~~~$Corrolary}
\newtheorem{coj}[theorem]{$~~~~$Conjecture}

\newtheorem{deff}[theorem]{$~~~~$Definition}

\def\nop{ }

\def\nor1{Normed$\{~2^{ \zzz \theta  \, )} ~$,$~\sqrt{~2^{ \zzz \theta  \, )}}~\}$}

\def\xor2{Normed$\{ ~\sqrt{~2^{ \zzz \theta  \, )}}~,~2~ \} $}
\def\glamb{\lambda}
\def\glamb{P}
\def\glamb{\theta}
\def\pag2{Page 2}

\def\glamb{\zeta}
\def\zzthe{\theta}

\def\zzz{~ \sharp ( ~ }

\def\f55{ \normalsize  \baselineskip = 1.8 \normalbaselineskip }

\def\f55{  \baselineskip = 1.1 \normalbaselineskip } 
\def\g55{  \baselineskip = 1.0 \normalbaselineskip } 
\def\s55{ \baselineskip = 1.0 \normalbaselineskip } 

\def\f55{  \baselineskip = 0.7 \normalbaselineskip }

\newcommand{\thx}[1]{Theorem \ref{#1}}
\newcommand{\cjx}[1]{Conjecture \ref{#1}}
\newcommand{\phx}[1]{Proposition \ref{#1}}
 \newcommand{\dfx}[1]{Definition \ref{#1}}
\newcommand{\lem}[1]{Lemma \ref{#1}}

\newcommand{\el}[1]{Line (\ref{#1})}

\newcommand{\eq}[1]{(\ref{#1})}

\newcommand{\underx}[1]{\overline{~ {#1} ~}}

\begin{document}

%\title{Rough Summary of March 2012 Research Before I Attended
%the AMS March 17-18 Conference}

%% \title{ On
%% How the Revival of a Diluted 
%% Version of
%% Hilbert's Consistency Program 
%% %is 
%% Should
%% Likely
%%  Be
%% % Plausible and
%% % Veru
%%  Germane
%% to Computer Science}

\title{On  How the Introducing of a 
 New $~\theta~$ Function Symbol
Into Arithmetic's Formalism Is 
Germane
to Devising Axiom Systems that Can 
Appreciate Fragments of Their Own
Hilbert Consistency}

\def\beq{\begin{equation}}
\def\enq{\end{equation}}

\def\bel{\begin{lemma}}
\def\enl{\end{lemma}}

\def\bec{\begin{corollary}}
\def\enc{\end{corollary}}

\def\bed{\begin{description}}
\def\ennd{\end{description}}
\def\bee{\begin{enumerate}}
\def\ene{\end{enumerate}}

\def\bxbxd{\begin{definition}}
\def\bxbxdd{\begin{definition}}
\def\eedd{\end{definition}}
\def\bxbxdr{\begin{definition} \rm}
\def\bel{\begin{lemma}}
\def\enl{\end{lemma}}
\def\ent{\end{theorem}}

\author{  Dan E.Willard\thanks{This research 
was partially supported
by the NSF Grant CCR  0956495.}}
%Email = dew@cs.albany.edu.}}
%\newline
%Email = dan.willard.albany@gmail.com}}

%%%\date{Copyright 2012 by Dan E. Willard}

\date{State University of New York at Albany}

\maketitle

\setcounter{page}{0}
\thispagestyle{empty}

\normalsize

\baselineskip = 1.3\normalbaselineskip

\normalsize

\baselineskip = 1.0 \normalbaselineskip 
\def\bbint{\large \baselineskip = 1.6 \normalbaselineskip } 
\def\bbint{\large \baselineskip = 1.6 \normalbaselineskip }
\def\bbint{\normalsize \baselineskip = 1.3 \normalbaselineskip }

%%\baselineskip = 1.0 \normalbaselineskip 
%%\def\bbint{\large \baselineskip = 1.6 \normalbaselineskip } 
%%\def\bbint{\large \baselineskip = 1.6 \normalbaselineskip }
%%\def\bbint{\normalsize \baselineskip = 1.3 \normalbaselineskip }

\def\bbint{\normalsize \baselineskip = 1.27 \normalbaselineskip }

\def\bbint{\large \baselineskip = 2.0 \normalbaselineskip }

\def\bbint{\normalsize \baselineskip = 1.25 \normalbaselineskip }
\def\bbina{\normalsize \baselineskip = 1.24 \normalbaselineskip }

\def\bbint{\large \baselineskip = 2.0 \normalbaselineskip }

\def\bbing{ }
\def\bbins{ }
\def\bbinm{ }

\def\bbint{\normalsize \baselineskip = 1.95 \normalbaselineskip }

\def\bbing{ }
\def\bbins{ }
\def\bbinm{ }

\def\bbint{\large \baselineskip = 2.3 \normalbaselineskip } 
\def\bbing{ }
\def\bbins{ }
\def\bbinm{ }

\def\bbint{\normalsize \baselineskip = 1.7 \normalbaselineskip } 

\def\bbint{\large \baselineskip = 2.3 \normalbaselineskip } 
\def\bbinm{ \baselineskip = 1.18 \normalbaselineskip }

\def\bbint{\large \baselineskip = 2.0 \normalbaselineskip } 
\def\bbing{ }
\def\bbins{ }
\def\bbinm{ }
\def\bbinr{ }

\def\bbint{\normalsize \baselineskip = 1.25 \normalbaselineskip }
\def\bbina{\normalsize \baselineskip = 1.24 \normalbaselineskip }
\def\bbinr{ \baselineskip = 1.3 \normalbaselineskip }
\def\bbing{ \baselineskip = 1.28 \normalbaselineskip }
\def\bbins{ \baselineskip = 1.21 \normalbaselineskip }
\def\bbinm{  }

\def\ftl{ \baselineskip = 1.5 \normalbaselineskip }

\bbint

\parskip 5 pt

\noindent

\small

\baselineskip = 1.14 \normalbaselineskip

\parskip 5pt

\baselineskip = 1.2 \normalbaselineskip 

%\setcounter{page}{0}

%%%%%%{\small

\large
\normalsize

 \baselineskip = 1.2 \normalbaselineskip 

%mmmmmmmmmmmmm

\begin{abstract}
\Large
\baselineskip = 1.65 \normalbaselineskip

A new $\theta$ function primitive is proposed that almost achieves the
combined efficiency of the addition, multiplication and successor growth
operations.
%  We conjecture this
This
 $\theta$ function symbol enables the
constructing of an ``IQFS(PA+)'' axiom system that can corroborate a
fragmentary definition of its own Hilbert consistency, while it will
simultaneously verify isomorphic counterparts of all Peano Arithmetic's
$\Pi_1$ theorems.
Many propositions and intermediate results are also
established.
Only one  intermediate result, which most readers
will intuit should be true, does remain formally unproven.
%  Only one very probable intermediate result is technically left unproven.

\end{abstract}

%% These considerations lead us to conjecture that new boundary-case
%% evasions of the Second Incompleteness Theorem exist, where an axiom system
%%  can
%% have a limited-but-real appreciation of a fragment of its own Hilbert
%% consistency. Apart from this conjecture, our new $\theta$ function primitive
%% is interesting in its own right.

%%% \leq{}
\bigskip
\bigskip


\bigskip


\bigskip

\small
{\bf Comment:$~$ } All the theorems and propositons are
the same in this Version 5 as in Version 4.
The difference is that the writing style  is now significantly
more polished.

%%{\bf Comment:$~$ } This version 4 consists
%%of 35  pages (plus a over page)  It is a
%%more 
%% refined result than our earlier drafts.

% 
% 
% 
% {\bf Mathematics Subject Classification:}
% 03B52; 03F25; 03F45; 03H13 
% 
% 
% 
% \bigskip
% \bigskip

% {\bf Please Cite this Paper as:}
% {\rm http://arxiv.org/abs/1108.6330}, 
%  appearing in Cornell Archives 

\def\ww22{\normalsize \baselineskip = 1.21\normalbaselineskip \parskip 4 pt}
\def\bb22{\normalsize \baselineskip = 1.19\normalbaselineskip \parskip 4 pt}
\def\zz22z{\normalsize \baselineskip = 1.19 \normalbaselineskip \parskip 3 pt}
\def\xx22{\normalsize \baselineskip = 1.17\normalbaselineskip \parskip 4 pt}
\def\vx22s{\normalsize \baselineskip = 1.16 \normalbaselineskip \parskip 3 pt} 
\def\vv22{\normalsize \baselineskip = 1.17 \normalbaselineskip \parskip 3 pt} 
\def\aa22{\normalsize \baselineskip = 1.15 \normalbaselineskip \parskip 3 pt} 
\def\g55{  \baselineskip = 1.0 \normalbaselineskip } 
\def\s55{ \baselineskip = 1.0 \normalbaselineskip } 
\def\sm55{ \baselineskip = 0.9 \normalbaselineskip }

\vspace*{- 1.0 em}

\def\waw11{\normalsize \baselineskip = 1.72\normalbaselineskip}
\def\waw11{\normalsize \baselineskip = 1.12\normalbaselineskip}
\def\waw11{\normalsize \baselineskip = 1.85\normalbaselineskip}

\def\waw11{\normalsize \baselineskip = 1.45\normalbaselineskip}

\def\waw11{\normalsize \baselineskip = 1.7\normalbaselineskip}

\def\waw11{\normalsize \baselineskip = 1.12\normalbaselineskip}

\def\g55{  \baselineskip = 1.50 \normalbaselineskip } 
\def\s55{ \baselineskip = 1.50 \normalbaselineskip } 
\def\sm55{ \baselineskip = 1.5 \normalbaselineskip }

\def\g55{  \baselineskip = 1.50 \normalbaselineskip } 
\def\s55{ \baselineskip = 1.50 \normalbaselineskip } 
\def\sm55{ \baselineskip = 0.9 \normalbaselineskip }

\def\aa22{\normalsize  \waw11 \parskip 6 pt} 
\def\bb22{\normalsize  \waw11 \parskip 5 pt}
\def\ww22{\normalsize \waw11 \parskip 4 pt}
\def\vv22{\normalsize  \waw11 \parskip 3 pt} 
\def\tt22{\normalsize  \waw11 \parskip 2 pt} 

\def\g55{  \baselineskip = 1.0 \normalbaselineskip } 
\def\b55{  \baselineskip = 1.0 \normalbaselineskip } 
\def\s55{ \baselineskip = 1.0 \normalbaselineskip } 
\def\sm55{ \baselineskip = 0.9 \normalbaselineskip }

\def\mal{ \bf  }
\def\nal{\mathcal}

\def\cvrew{ \baselineskip = 1.6 \normalbaselineskip \parskip 3pt }

\def\ttt2c{ }
\def\tttc{ }

\def\tttc{\tiny \baselineskip = 0.8 \normalbaselineskip  \parskip 0pt }
\def\ttt2c{\tiny \baselineskip = 0.7 \normalbaselineskip  \parskip 0pt }
\def\tttc{ \baselineskip = 2.1 \normalbaselineskip  \parskip 5pt }
\def\ttt2c{ \baselineskip = 2.1 \normalbaselineskip  \parskip 5pt }

\def\tttc{ \baselineskip = 1.15 \normalbaselineskip  \parskip 5pt }
\def\ttt2c{ \baselineskip = 1.15 \normalbaselineskip  \parskip 5pt }

\def\tttc{ \baselineskip = 1.12 \normalbaselineskip  \parskip 4pt }
\def\ttt2c{ \baselineskip = 1.12 \normalbaselineskip  \parskip 4pt }

\def\tttc{ \baselineskip = 1.14 \normalbaselineskip  \parskip 3pt }
\def\ttt2c{ \baselineskip = 1.14 \normalbaselineskip  \parskip 4pt }

\def\cvt{ \baselineskip = 0.98 \normalbaselineskip }
\def\cv9{ \baselineskip = 0.99 \normalbaselineskip }
\def\cvs{ \baselineskip = 1.0 \normalbaselineskip }
\def\cvl{ \baselineskip = 1.0 \normalbaselineskip }
\def\cvh{ \baselineskip = 1.03 \normalbaselineskip }
\def\cvg{ \baselineskip = 1.00 \normalbaselineskip }

\def\cvt{ \baselineskip = 1.6 \normalbaselineskip }
\def\cv9{ \baselineskip = 1.6 \normalbaselineskip }
\def\cvs{ \baselineskip = 1.6 \normalbaselineskip }
\def\cvl{ \baselineskip = 1.6 \normalbaselineskip }
\def\cvh{ \baselineskip = 1.6 \normalbaselineskip }
\def\cvg{ \baselineskip = 1.6 \normalbaselineskip }
\def\cvb{ \baselineskip = 1.6 \normalbaselineskip }
\def\cvnew{ \baselineskip = 1.6 \normalbaselineskip }
\def\cvmew{ \baselineskip = 1.6 \normalbaselineskip }
\def\cvwew{ \baselineskip = 1.6 \normalbaselineskip \parskip 5pt }
\def\cvrew{ \baselineskip = 1.6 \normalbaselineskip \parskip 3pt }

\def\cvt{ \baselineskip = 1.22 \normalbaselineskip }
\def\cv9{ \baselineskip = 1.22 \normalbaselineskip }
\def\cvs{ \baselineskip = 1.22 \normalbaselineskip }
\def\cvl{ \baselineskip = 1.22 \normalbaselineskip }
\def\cvh{ \baselineskip = 1.22 \normalbaselineskip }
\def\cvg{ \baselineskip = 1.22 \normalbaselineskip }
\def\cvb{ \baselineskip = 1.22 \normalbaselineskip }
\def\cvnew{ \baselineskip = 1.4 \normalbaselineskip }
\def\cvmew{ \baselineskip = 1.35 \normalbaselineskip }
\def\cvwew{ \baselineskip = 1.4 \normalbaselineskip \parskip 5pt }
\def\cvrew{ \baselineskip = 1.22 \normalbaselineskip \parskip 3pt }

\def\cvt{ }
\def\cv9{ }
\def\cvs{ }
\def\cvl{ }
\def\cvh{ }
\def\cvg{ }
\def\cvb{ }
\def\cvnew{ } 
\def\cvmew{ }
\def\cvwew{ }
\def\cvrew{ }

\def\fend{ 

\medskip -------------------------------------------------------}

\def\g55{  \baselineskip = 1.0 \normalbaselineskip } 
\def\s55{ \baselineskip = 1.0 \normalbaselineskip } 
\def\sm55{ \baselineskip = 1.0 \normalbaselineskip } 
\def\h55{  \baselineskip = 1.08 \normalbaselineskip } 
\def\b55{  \baselineskip = 1.1 \normalbaselineskip } 

\normalsize

\baselineskip = 1.85 \normalbaselineskip

%% Sleepy  

%\cvlpm %% Sleepy  
%\cvnew

%\small

%\parskip 0p

\parskip 2pt

\vspace*{- 1.0 em}

% \newpage

%\large

%\setcounter{page}{0}
\baselineskip = 1.04 \normalbaselineskip 
\parskip 2pt

\baselineskip = 0.96 \normalbaselineskip 
%\baselineskip = 0.90 \normalbaselineskip 

%\parskip 1pt
% 
\baselineskip = 2.16 \normalbaselineskip 
\baselineskip = 2.3 \normalbaselineskip 

\baselineskip = 0.95 \normalbaselineskip 
%\baselineskip = 0.95 \normalbaselineskip 
\baselineskip = 0.88 \normalbaselineskip 
\parskip 0pt

\noindent

% 
% 
% NNEW COMMENT
% 
% 
% The pdf version of this draft is verbatim identical to August's Version 3.
% The prior draft's abstract was incorrectly broadcast by Arxiv on the 
% Internet, after I pressed a wrong computer button. Thus, 
% Version 4 was issued.

\newpage

\def\gvs{ \normalsize \baselineskip = 1.4 \normalbaselineskip  \parskip    5pt}
\def\gvs{ \normalsize \baselineskip = 1.44 \normalbaselineskip  \parskip    5pt}
\def\gvs{ \large \baselineskip = 1.44 \normalbaselineskip  \parskip    5pt}
\def\gvs{ \normalsize \baselineskip = 1.44 \normalbaselineskip  \parskip    5pt}\def\gvs{ \normalsize \baselineskip = 1.74 \normalbaselineskip  \parskip    5pt}
\def\gvs{ \normalsize \baselineskip = 1.44 \normalbaselineskip  \parskip 5pt}

\def\gvs{   \baselineskip = 1.74 \normalbaselineskip  \parskip    5pt}

\def\gvs{ \normalsize \baselineskip = 1.44 \normalbaselineskip  \parskip 5pt}
\def\gvs{ \large \baselineskip = 2.0 \normalbaselineskip  \parskip 5pt}
\def\gvs{ \Large \baselineskip = 2.0 \normalbaselineskip  \parskip 5pt}
\def\gvs{ \normalsize \baselineskip = 2.44 \normalbaselineskip  \parskip 5pt}
\def\gvs{ \normalsize \baselineskip = 2.04 \normalbaselineskip  \parskip 5pt}
\def\gvs{ \normalsize \baselineskip = 2.64 \normalbaselineskip  \parskip 5pt}
\def\gvs{ \Large \baselineskip = 1.6 \normalbaselineskip  \parskip 5pt}

\gvs

\footnotesize

\def\gvs{ }

\normalsize \baselineskip = 0.98 \normalbaselineskip
\normalsize \baselineskip = 1.0 \normalbaselineskip
\normalsize \baselineskip = 1.01 \normalbaselineskip

\def\gvs{ \normalsize \baselineskip = 1.25 \normalbaselineskip  \parskip 4pt}

\def\gvs{ \Large \baselineskip = 1.6  \normalbaselineskip  \parskip 6pt}
\def\gvs{ \normalsize \baselineskip = 1.6  \normalbaselineskip  \parskip 6pt}
\def\gvs{ \large \baselineskip = 1.6  \normalbaselineskip  \parskip 6pt}

\def\gvs{ \normalsize \baselineskip = 1.227 \normalbaselineskip  \parskip 3pt}
\def\gvs{ \large \baselineskip = 1.8  \normalbaselineskip  \parskip 6pt}

\def\gvs{ \normalsize \baselineskip = 1.5 \normalbaselineskip  \parskip 3pt}

\def\gvs{ \large \baselineskip = 2.1  \normalbaselineskip  \parskip 6pt}

\def\gvs{ \normalsize \baselineskip = 2.1  \normalbaselineskip  \parskip 6pt}

%%%old 

 \def\gvs{ \normalsize \baselineskip = 1.227 \normalbaselineskip  \parskip 3pt}

 \def\gvs{ \large  \baselineskip = 1.6 \normalbaselineskip  \parskip 5pt}
%% march 31

\def\gvs{ \Large  \baselineskip = 1.8 \normalbaselineskip  \parskip 5pt}
\def\gvs{ \LARGE  \baselineskip = 1.8 \normalbaselineskip  \parskip 5pt}
\def\gvs{ \normalsize  \baselineskip = 2.0 \normalbaselineskip  \parskip 5pt}

\def\gvs{ \Large  \baselineskip = 2.0 \normalbaselineskip  \parskip 5pt}

\def\gvs{ \large  \baselineskip = 2.2 \normalbaselineskip  \parskip 5pt}

\def\gvs{ \normalsize \baselineskip = 2.4  \normalbaselineskip  \parskip 6pt}

\def\gvs{ \normalsize \baselineskip = 2.6  \normalbaselineskip  \parskip 6pt}
\def\gvs{ \normalsize \baselineskip = 2.2  \normalbaselineskip  \parskip 6pt}
\def\gvs{ \normalsize \baselineskip = 1.8  \normalbaselineskip  \parskip 5pt}

\def\sgvs{ \small \baselineskip = 1.33  \normalbaselineskip  \parskip 1pt}
\def\tttc{ }
%\baselineskip = 1.14 \normalbaselineskip  \parskip 4pt }
\def\ttt2c{ }
%\baselineskip = 1.14 \normalbaselineskip  \parskip 4pt }

\def\gv2{ \normalsize \baselineskip = 1.30  \normalbaselineskip  \parskip 3pt}

\def\gvs{ }

\def\gvs{ \normalsize \baselineskip = 2.1 \normalbaselineskip  \parskip 7pt}
\def\gvs{ \normalsize \baselineskip = 1.8 \normalbaselineskip  \parskip    7pt}

 \def\gvs{ \large \baselineskip = 1.7  \normalbaselineskip  \parskip 9pt}
\def\gvs{ \normalsize \baselineskip = 2.0  \normalbaselineskip  \parskip 9pt}

\def\gv2{ \normalsize \baselineskip = 1.30  \normalbaselineskip  \parskip 3pt}

\def\gvs{ \large \baselineskip = 1.7  \normalbaselineskip  \parskip 5pt}

\def\gvs{ \normalsize \baselineskip = 2.0  \normalbaselineskip  \parskip 8pt}
\def\gvs{ \large \baselineskip = 2.0  \normalbaselineskip  \parskip 8pt}

%fffff

%fffff
\def\gvs{ \normalsize \baselineskip = 1.3 \normalbaselineskip  \parskip   5pt}

\def\gvx{ \normalsize \baselineskip = 1.23 \normalbaselineskip  \parskip    3pt}

\section{Introduction}
% \section{Scientific Notes of Dan Willard Notarized on Nov 22,2016}
%%%%%%%%%% 1111111111111111}
\label{ss1}

\gvs
\tttc

Two historic results
were established
by
%%% in
 G\"{o}del's
millennial
%centennial
%  1931 
 paper \cite{Go31}.  
The First Incompleteness
Theorem 
showed  
no decision procedure
exists for identifying
all the true statements of Arithmetic.
%G\"{o}del's
The Theorem XI in \cite{Go31},
later known as the ``Second  Incompleteness
Theorem'',
%,
%appearing in G\"{o}del's millennial paper \cite{Go31}. 
demonstrated
% that
no extension
of
% axiom systems,
% roughly corresponding to
the
 Russell-Whitehead Principia Mathematicae formalism
% $\, P \,$
can
% could
 verify
its own consistency.
 G\"{o}del's
two
observations
% Theorem XI 
%sss%
 were  
%are
historic
% mainly 
 because they
 demonstrated, unequivocally,
that the 
initial 
objectives of Hilbert's Consistency
Program were 
much
too far-reaching.
Thus at best,  only a
% very  
sharply curtailed
form of Hilbert's goals 
would be
% was 
plausible.

These observations were
% This fact was 
%%further 
reinforced by a new version of the
Second Incompleteness Theorem, due to
the combined work of Pudl\'{a}k and Solovay \cite{Pu85,So94},
enhanced by some added techniques of Nelson and Wilkie-Paris
\cite{Ne86,WP87}.
 It 
% They collectively
established the prohibition
(stated in
% (in
 \textsection \ref{ss3}'s Theorem ++) that
no conventional
arithmetic formalism
can verify its own 
``Hilbert-style''
consistency, when it merely
recognizes Successor as a 
formally
total
functional operation.

%plausibly feasible
%from Theorem XI's perspective.

Within such curtailed limits, we have published since 1993
a series of articles 
\cite{ww93}-\cite{ww14},
outlining generalizations of the Second Incompleteness
Theorem and
its 
sometimes-feasible
% plausible 
 boundary-case exceptions.
%  that were formally feasible.
Pavel Pudl\'{a}k
examined
%in great detail,
% the preprints 
an early preprint
 of  
our
article \cite{wwapal}
and suggested 
\cite{Pupriv} 
%that 
we
consider  attempting to hybridize 
its formalism with some of
Ajtai's observations about 
the Pigeon Hole effects
\cite{Aj94}.
On a more informal basis during a lunch at a 1997 conference,
%meeting, 
Sam Buss \cite{Bu97} 
conveyed to us an approximately similar
suggestion.  
The Section 6 of 
  \cite{wwapal}
did subsequently
formalize one type of response to 
these insightful observations.
Our new results in this current paper will examine
a much more ambitious approach to \cite{wwapal}'s 
particular
topic,
that is germane to 
an arithmetic appreciating
some very
{\it delicately defined}
fragments of
its Hilbert-styled formalizations 
of its self
consistency.

During our discussion,
 $Add(x,y,z)$ and    $Mult(x,y,z)$ 
 will
denote 
two 
3-way predicate symbols
specifying
that
$x+y=z$ and
$x*y=z$,
which are known by our formalism to 
satisfy $\Pi_1$ encodings for axioms
formalizing the associative, commutative, identity and
distributive principles.
Let us say 
%that 
an
axiom 
system
%basis
 $\, \alpha \,$
{\bf recognizes}
successor, 
 addition  and multiplication
as {\bf Total Functions} iff 
%%%%%%%%% $\, \alpha \,$ 
it
can 
additionally
prove
\eq{totxtefs} - \eq{totxtefm}
as theorems.

 {
\sgvs
\vspace*{- 0.6 em}

{
%\small
\beq 
\label{totxtefs}
\forall x ~ \exists z ~~~Add(x,1,z)~~
\enq
\vspace*{- 1.7 em}
\beq 
\label{totxtefa}
\forall x ~\forall y~ \exists z ~~~Add(x,y,z)~~
\enq
\vspace*{- 1.7 em}
\beq 
\label{totxtefm}
\forall x ~\forall y ~\exists z ~~~Mult(x,y,z)~
\enq }
}

\vspace*{- 1.4 em}
\noindent
Our 
studied
axiomatizations for integer arithmetics 
%in the current article 
will differ from 
their conventional counterparts
{\it by neither 
possessing} an ability to prove any of the totality statements
\eq{totxtefs} - \eq{totxtefm} above, 
{\it nor
containing} function symbols for
formalizing
the traditional addition, multiplication and successor 
function
operations.
%%% nor recognizing the validity of the identities
% 
% which state that these primitive
% are ``total functions''. 
% 
Instead, we will rely upon an entirely
new
different type of primitive function,
% symbol, 
called the ``$~\theta~$''
operator, to construct the endless sequence of integers
$~3,4,5, ~...~$ from
merely the
three  initial starting constants of 0, 1 and 2.
(The exact number of logical symbols for encoding a term $T_n~$,
that represents an integer $\, n \geq 3 \,$, will satisfy
 either an O( $Log^3n$ ) or  O( $Log \, n$ ) upper bound
under our 
Propositions \ref{th-3.3} and
\ref{th-7.1}, 
depending on whether a tree, as opposed to 
a Dag-oriented methodology,
is used to encode $T_n$ .)

Assuming only
% that 
the very modest assumptions
of
% that
our \cjx{con6.6} 
are true$\,$\footnote{\cjx{con6.6} 
is the 
% formally 
%%%%%% unproven
% intermediate 
result, mentioned
% noted 
in the abstract, 
which shall be 
demonstrated to be  almost certainly  true
by the
appendix attached to \textsection\ref{nn6}. } ,
%% 
%%which we are quite
%%certain most readers
%%% will 
%%shall
%%agree is true.}  ,
%%
% our 
Theorem \ref{thmain} 
will imply
any consistent arithmetic axiom system $A$ (including Peano
Arithmetic
{\it itself} ) can be mapped onto a consistent axiom system
``IQFS$(A)$'', where  IQFS$(A)$ can {\it simultaneously} :
\bee
% \small \baselineskip = 1.32  \normalbaselineskip  \parskip 2pt
\item  
corroborate its own Hilbert consistency,
\item 
confirm the validity of all $A$'s $\Pi_1$ theorems in 
a slightly revised
% modified
language that uses the preceeding  
 $Add(x,y,z)$ and    $Mult(x,y,z)$  3-way predicate symbols, and
\item 
use only the three starting constants of 0, 1 and 2, 
to construct any term $T_n$ for representing any natural
number $n \geq 3$, while satisfying
Propositions \ref{th-3.3} and
\ref{th-7.1}'s
specified
% previously mentioned
 O($ \, Log^3n\,$) and  O($\,Log \, n$)
  bounds for the length of $T_n$'s construction.
\ene
Items
%The combination of
% Items 
1-3 are potentially
significant
% noteworthy
% highly  significant
 because \cite{wwapal}'s 
% less advanced
earlier
ISCE formalism 
could  achieve the first two objectives, only when it
had   used
an
% unfortunately undesirable 
unfortunately
infinite number of 
% distinct
constant symbols. 

Theorem \ref{thmain}'s
result should be 
viewed,
cautiously,
  {\it as  only a boundary-case
evasion} of the Second Incompleteness Theorem.
%% {\it even if our employed \cjx{con6.6}} is correct.
This is
 because
its  ``IQFS(A)'' formalism will  be too weak to
 treat
even
 Successor as a total function. Yet assuming our
\cjx{con6.6} is correct
(as the Appendix will
% explain 
demonstrate
it is nearly 100 \%
% likely 
certain
to be), the semantics of
IQFS shall
partially answer an open question, 
that was
raised by Harvey Friedman,
at the end of the fourth lecture within his 5-part
  year-2014 You Tube 
%lecture 
%presentation
series \cite{Fr14}.

It concerned
whether some 
% unusual 
%types of 
% non-trivial
exotic
 boundary case
exceptions to the Second Incompleteness Theorem might exist. 
The intuition behind 
IQFS(A)'s
%  new  results
improved behavior 
will be 
%is
 that
Proposition \ref{th-3.3}'s
% newly 
% proposed
$~\theta~$ primitive
will  capture much of
the growth properties
of the 
%%%%%traditional 
addition, multiplication and successor operations,
{\it without employing
those of their special functionalities,} 
that are known to particularly
% directly 
support 
\cite{Pu85,So94}'s 
broad-scale
%unique
generalizations
 of
the Second Incompleteness Theorem.

% ccc

\medskip

 The present
 article,
fortunately, 
%has been carefully composed
%so that it 
can be
read 
 without 
% first
examining any of our previous
papers. If a reader
does wish to skim
one of our earlier
articles,
% before the current  article,
we recommend 
Sections 3 
\& 4
 of \cite{wwapal}.
(The reading of 
\cite{wwapal} is
% actually, 
%technically,unnecessary
optional
because a later chapter will
% adequately
 summarize its content in more than
adequate detail.)

\bigskip

%% \medskip
%% 
%% We suspect
%% % that 
%% %% most experts about Bounded Arithmetic will conclude, 
%% that most researchers into proof theory
%% will agree,
%% by the end of this paper, that our 
%% \cjx{con6.6} is likely to be correct. 
%% Therefore, some types of
%% significant
%%  boundary-case
%% evasions of the 
%% Second Incompleteness Theorem 
%% will seemingly arise.

%are feasible. 

\section{Returning to the 1931-1939 Period}
\label{ss2}
 G\"{o}del's 
Second Incompleteness Theorem 
was published in two 
%quite different 
% forms
stages
 during the
1931-1939 period.
Its initial 1931
% variant, 
announcement,
formalized by Theorem XI
in  G\"{o}del's millineal paper \cite{Go31},
% 
% 
% Its Theorem XI,
% later known as the ``Second  Incompleteness
% Theorem'',
% %,
% %appearing in G\"{o}del's millennial paper \cite{Go31}. 
% 
demonstrated
that
no extension
of
% axiom systems,
% roughly corresponding to
the
 Russell-Whitehead Principia Mathematicae formalism
% $\, P \,$
could
% could
% verify
corroborate
its own consistency. 
The widely quoted more general
result, that 
every consistent r.e. extension
 of Peano Arithmetic must
be unable to 
%ss% prove a theorem affirming its
corroborate its
own consistency, 
was, technically,
 first
published 
in the 1939 edition of 
the Hilbert-Bernays
textbook \cite{HB39}. The latter has been considered
to be the 
first
definitive\footnote{\label{bolc}    Boolos \cite{bool}
states there has been some debate among scholars
whether  Peano Arithmetic's 
precise
generalization
of the Second Incompleteness Theorem should be fully credited
to G\"{o}del's seminal paper \cite{Go31}, as opposed
to having been implicit from it.
In any case, the Hilbert-Bernays
textbook \cite{HB39}
 vented this generalization 
% of it 
in its 1939 second edition, and 
von Neumann (unpublished) is also known
\cite{Da97,Go5,Yo5}
 to have made
similar observations.}
 demonstration 
of the broad reach of
the Second Incompleteness Effect. 
Its   
%generalized 
formalism 
established, beyond any reasonable doubt, that any
kind of axiom system
% type of formalism 
%possessing a conventional knowledge of 
corroborating
its own consistency
{\it must rely upon a 
foundational structure
%  fundamentally 
different from Peano Arithmetic.} 
(This is because the 
Hilbert-Bernays 
textbook formalized the forerunner of
what is now  known as the 
Hilbert-Bernays Derivability Conditions \cite{HB39,HP91,Lo55,Mend},
as a mechanism for
% foreseeing 
envisioning
the 
% astonishing
broad generality of the
Second Incompleteness Effect.)

It is, thus, fascinating that Hilbert,
as a co-author of 
an important 
%very 
% historic
generalization of the Incompleteness Theorem,
never withdrew the
% chose to  never fully withdraw his 
%1926 
justification
 \cite{Hil26}
for his consistency program:

\begin{quote}
\small
\baselineskip = 1.2 \normalbaselineskip
\ttt2c 
$*~$
{\it ``
Let us admit that the situation in which we presently
find ourselves with respect to paradoxes is in the long
run intolerable. Just think: in mathematics, this paragon of
reliability and truth, the very notions and inferences,
as everyone learns, teaches, and uses them, lead to absurdities.
And 
where 
else 
would 
reliability and truth be found 
if  even mathematical thinking fails?''}
\end{quote}
Indeed, 
%Instead,
Hilbert 
always insisted that some
special new formalism would at
least partially vindicate the
prior
%initial 
 goals of his
consistency program.
He thus arranged to 
have  its often-quoted 
motto
({\it ``Wir m\"{u}ssen wissen---Wir werden wissen''} )  
%of his 
%nevertheless, 
engraved
on his tombstone \footnote{ English translation: {\it
``We must know, We will know.''}}.

%Moreover, it 

It is 
known
\cite{Da97,Go5,Yo5}
that G\"{o}del
was
also
doubtful about the generality of the Second Incompleteness
Theorem for at least two years after its publication.
He thus inserted the following
%cautious 
% very famous
historically noteworthy
caveat into
his
famous 1931 
millennial 
paper  \cite{Go31}:
% whose closing section 
%%% 
%%% One of the closing paragraphs of
%%%  \cite{Go31} 
%%% thus
%%% included
%
%
%%% contained the following cautious disclaimer:
%caveat:
\begin{quote}
\small
\baselineskip = 1.2 \normalbaselineskip
\ttt2c 
$~**~~$ 
{\it ``It must be expressly noted that
Theorem XI''}
(e.g. the Second Incompleteness Theorem) 
{\it ``represents no contradiction of the formalistic
standpoint of Hilbert. For this standpoint
presupposes only the existence of a consistency
proof by finite means, and {\it there might
conceivably be finite proofs} which cannot
be stated in P or in ... ''}
\end{quote}

The 
% above 1931
statement $**$ has 
had
%been subject to 
numerous
%many
different
interpretations
\footnote{
Some 
scholars
have interpreted
$\,**\,$
as
%as, possibly,'
anticipating attempts
to confirm Peano Arithmetic's consistency,
via 
either
Gentzen's formalism or 
 G\"{o}del's Dialetica interpretation.
Other scholars have viewed $**$ 
as having more ambitious goals in 1931 (that is, goals which  G\"{o}del
 essentially  formally withdrew later).
}.
All
 G\"{o}del's 
biographers
\cite{Da97,Go5,Yo5}
%%%have
have noted
his
% that G\"{o}del's
 intention (prior to 1930)
was
to
establish
%achieve
Hilbert's proposed objectives, before 
he proved
%proving
% G\"{o}del proved
a result 
% 
% however,
% %%%%%his
% G\"{o}del
% did  originally 
% seek 
% % goal was 
% to
% establish
% %achieve
% Hilbert's proposed objectives before 
% proving
% % G\"{o}del proved
% a result 
% 
leading
%that led 
in the opposite direction.
Yourgrau \cite{Yo5}
records
%furthermore,
 how
von Neumann 
surprisingly
%did
{\it ``argued 
against G\"{o}del 
himself''}
in the early 1930's,
 about the definitive 
%%%  achievement of a'
 termination of Hilbert's
consistency program,  
which
{\it ``for several years''} after \cite{Go31}'s publication,
G\"{o}del 
{\it ``was cautious not to prejudge''}.
It is known 
 G\"{o}del 
began to 
more fully
endorse 
the Second Incompleteness
Theorem
during a 1933
%% Vienna
lecture  \cite{Go33},
and he 
% told biographers he
completely embraced it
after learning about Turing's work
\cite{Tur36}.

Our research
in \cite{ww93}-\cite{ww14}
%has been 
is
related to issues 
%analogous
 similar 
to those
%that were 
raised by Hilbert and 
G\"{o}del
in
 statements  $*$ and $**$.
This is because it is
% counter-intuitive and 
awkward to
%presume that 
explain how
human beings can maintain the
psychological drive and 
needed energy-desire to cogitate, without
being stimulated by an instinctive faith in their own
consistency (under a definition of 
% formal  consistency
such
% this concept 
that is suitably 
gentle and
% delicate
soft
 to 
be consistent with the 
Incompleteness 
Effect).

%Theorem's requirements).

% preclude a violation of the 
% restrictions imposed by the Incompleteness Theorem).

Accordingly, our research in
\cite{ww93}-\cite{ww14} 
has explored both generalizations and
boundary-case exceptions for the Incompleteness Effect, so as
to determine 
{\it exactly
%what type} 
which}
boundary-case evasions are 
permissible.
% permitted.
Our prior research in \cite{ww93}-\cite{ww14} 
% had 
used mostly cut-free forms of deduction to
evade the 
restrictions imposed by the
Second Incompleteness Effect. The current article will instead
focus on the more pristine Hilbert-Frege methods of deduction.
%xxx
Assuming the tempting assumptions of our
\cjx{con6.6} are correct (as we do), 
the present article will show
how a new type of 
 $~\theta~$  function symbol, 
employing an 
% unusual 
{\it ``indeterminate-styled''} definition of growth
(formalized in \textsection \ref{ss4}),
will usher in a surprising evasion
of the Second Incompleteness
Effect.

%n% 
%n% They
%n% % are likely
%n% will be shown to
%n% support an evasion of the Second Incompleteness
%n% Effect when our axiom systems replace the traditional  
%n% growth properties of the addition, multiplication and successor 
%n% function symbols with a new $~\theta~$ primitive
%n% that has an quite novel  ``indeterminate-styled'' definition of growth
%n% (as will be formalized in \textsection \ref{ss4}.)
%n% 
%n% The motivation for this replacement will be
%n% explained during the next 
%n% two
%n% sections of this article.

More details about $\theta$ will be explained
in the next two sections. 
Essentially, 
%it will be 
  $\theta$ is
needed
% essentially
because
a
 version of the Second Incompleteness Theorem,
due to the
combined work of Pudl\'{a}k, Solovay, Nelson and Wilkie-Paris
\cite{Ne86,Pu85,So94,WP87},
%will show 
demonstrates
that 
if an axiom system $~\alpha~$
proves
merely
 any one 
of \eq{totxtefs} - \eq{totxtefm}'s totality statements,
then it will be incapable of confirming its own consistency
under a Hilbert-style deductive method.

%%%%%     \smallskip

Our   Theorem \ref{thmain} 
 will 
suggest  it is possible to obtain a
%ss% {\it part-way 5-10 \% positive} 
{\it partially positive}
interpretation
for what
Hilbert and G\"{o}del 
were
seeking 
% a Consistency Program
% to establish
%%  seeking to accomplish
% contemplating 
in their 
statements
$*$ and $**$, 
in
% within
 a context where 
it is  known 
% that 
the 
Second Incompleteness Effect,
clearly,
 precludes a full achievement of
their
earlier
% Hilbert's
objectives.

\section{Motivation for Research and
its Broad Perspective}
% and Background Notation}
 %  222222}
\label{ss3}

It is helpful to employ a flexible vocabulary so
% that our
our 
results
will apply 
%research 
%to the 
%formalisms
to any of the
%% 
%% accessible to
%% some 
%% readers who are acquainted with only one of the 
%% 
textbook 
formalisms of
% settings outlined by 
%  
say 
Enderton,
Fitting, H\'{a}jek-Pudl\'{a}k,
or Mendelson
\cite{End,Fit,HP91,Mend}.
% 
% ,
% %% or
% %%  Papadimitriou 
% %% \cite{End,Fit,HP91,Mend,Papa},
% %% 
% %% 
% %the widest
% % possible  
% %audience, 
% it is 
% helpful
% % useful to 
% use
% %employ 
% a
% % very 
% flexible
%  vocabulary.
% %%
% %%that 
% %%allows a reader to 
% %%%quickly 
% %%%translate results 
% %%traverse
% %%from
% %%one textbook to another.
% % Therefore, let us define a
Let us 
% thereby
call an
ordered pair $(\alpha,d)$ a
    {\bf ``Generalized Arithmetic''}
% therefore
iff its 
% first and second 
two components 
are 
% described
% formalized
defined 
as 
follows:
\bee
\item
The {\bf Axiom Basis} ``$~\alpha~$'' 
of an 
% arbitrary 
arithmetic
%shall be 
is
defined as
its set of 
{\it proper axioms} 
it employs.
% employed by the formalism  $(  \alpha  , d  )$.
\item
% An arithmetic's
The  {\bf Deductive Apparatus} ``$~d~$''
of an arithmetic 
is defined as
the 
{\it combination} of its formal rules for inference
and 
%its
the
 built-in
 logical axioms ``$~L_d~$''
%  (that are 
% implicitly 
that are 
used
% employed 
by these rules.
\ene

\begin{exx}
\label{ex-2.1}
%\label{ex-basis}
\rm
This notation 
allows one to
 conveniently separate  the logical axioms
$~L_d~,~$ associated
with  $(  \alpha  , d  )~$, from 
%ss% $\, \alpha \,$'s
the 
``basis axioms'' $\, \alpha \,$.
%basis axioms
It also allows one to  isolate
and compare
%  , conveniently, 
various
apparatus techniques,
%technique,
% employed in the exact formalisms 
including the
 $~d_E~$,  
 $~d_M~$, 
 $~d_H~$,  
and  $~d_F~$
methods
%that we will now define:
defined below:
%% 
%%  Three
%% examples of this are illustrated below,
%% in a context where
%% are the deductive apparatus machineries defined
%% in  Enderton's, Mendelson's and Fitting's textbooks
%% \cite{End,Fi96,Mend}.
%% 
\bed
\item[   i. ]
The  $~d_E~$ apparatus,
%  formalized
introduced
 in 
\textsection
 2.4 of  Enderton's textbook, 
% will
uses only  modus ponens
as a rule of inference.
The latter will be accompanied
by 
a
4-part 
system of
 logical axioms,
called $~L_{d_E}~$, $\,$ to endow
 $~d_E~$ 
with an
ability to  support
G\"{o}del's Completeness Theorem.
\item[  ii. ]
The  $~d_M~$ 
apparatus in
\textsection 2.3 
of  Mendelson's textbook
and  the $d_H$ 
 apparatus
in  \textsection 0.10
of the H\'{a}jek-Pudl\'{a}k's
 textbook
employ a more compressed set of logical axioms
than $\, d_E \,$,
in a context where
they use
two rules of inference
(modus ponens and generalization).
%% plus a smaller set of logical axioms, which Mendelson
%% has called A1-A5.
%% Also, the $d_H$ 
%%  apparatus 
%% on pages ???? 
%% of the 
%%  H\'{a}jek-Pudl\'{a}k textbook
%% uses a slightly different variation of a generalization. 
%% (In the end, 
In the end, 
% both
 $~d_M~$ 
and  $~d_H~$ 
prove the same set of theorems
as   $~d_E~$ with 
only minor and unimportant changes in
proof length.
\item[ iii. ]
The 
``semantic tableaux''
 $~d_F~$ 
apparatus in
Fitting's 
and Smullyan's 
textbooks 
\cite{Fit,Smul}
was
 the main focus of our 
investigations in \cite{ww93,ww1,wwlogos,ww5,ww6,ww14}.
It will be rarely used
in the current article,
however.
Unlike
 $~d_E~$,  $~d_M~$  and  $~d_H~$, it
employs no logical axioms.
It instead
 uses  more complicated rules of inference.
It has many applications in
%% 
%% This tableaux apparatus 
%% % and also Resolution, have  been 
%%  has 
%% % been found to  have many 
%% a wide array of
%% applications
%% for 
%% 
automated deduction,
%although it is 
but is
less efficient than
 $  d_E  $,  $  d_M  $  and  $  d_H  $   
in
% under
% extremal
worst-case 
environments.
% settings.
%circumstances.
\ennd
\end{exx}

\begin{dff}
\label{def-2.2}
\rm
The  term {\bf ``Hilbert-style''} deductive method shall
refer to any  deductive  apparatus $~d~$ which
employs a modus ponens rule and also satisfies
G\"{o}del's Completeness Theorem.
(Thus, each of the  $  d_E  $,  $  d_M  $  and  $  d_H  $  
are  examples of Hilbert-style deductive methodologies.)
\end{dff}

\begin{exx}
\label{ex-2.3}
%\label{ex-basis}
\rm
Some added notation is
 needed to
explain why
% help outline
% an important distinction between 
a Hilbert style
deductive apparatus, such as $\,d_E\,$,  $\,d_H\,$
 or $\,d_M\,$, should be distinguished from
 $d_F$'s
``tableaux'' apparatus.
Let
% the symbols
 $Add(x,y,z)$ and    $Mult(x,y,z)$
%once 
again 
% will
denote 
two 
3-way predicate symbols
specifying
that
$x+y=z$ and
$x*y=z$.
%%%%%%%%
Also, let 
%Let
us recall 
% that 
an
axiom basis
 ``$\, \alpha \,$''
is said to
{\bf recognize}
successor, 
 addition  and multiplication
as {\bf ``Total Functions''} iff 
%%%%%%%%% $\, \alpha \,$ 
% it includes
% Lines 
\eq{totdefxs}-\eq{totdefxm} are
among its 
%  set of 
derived theorems.

% {\small
{\vspace*{- 0.6 em}
{ 
\sgvs
%\small
\beq 
\label{totdefxs}
\forall x ~ \exists z ~~~Add(x,1,z)~~
\enq
\vspace*{- 1.7 em}
\beq 
\label{totdefxa}
\forall x ~\forall y~ \exists z ~~~Add(x,y,z)~~
\enq
\vspace*{- 1.7 em}
\beq 
\label{totdefxm}
\forall x ~\forall y ~\exists z ~~~Mult(x,y,z)~
\enq }
}

\vspace*{- 1.4 em}

\noindent
%Then an 
%%%In this context, an
An
``axiom basis''
$\alpha$
%will be
is
 called 
{\bf Type-M} iff it includes
\eq{totdefxs}-\eq{totdefxm}
% \ref{totdefxs}-\ref{totdefxm}
as theorems,  
{\bf Type-A} if it includes 
%only
\eq{totdefxs} and \eq{totdefxa} as theorems,
and {\bf Type-S} if it contains
only \eq{totdefxs} as a
 theorem. 
%Moreover, 
Also,
$\alpha$ 
% will be 
is
called 
{\bf Type-NS} iff it  can prove
none of these theorems.
A summary of the prior literature, using
this particular
% the preceding 
notation,
is provided below:
\smallskip
\bed
\item[   $~~~~$a.$~~$  ]
The
combined research of Pudl\'{a}k, Solovay, Nelson and Wilkie-Paris
\cite{Ne86,Pu85,So94,WP87},
as 
is
formalized by statement $\, ++ \,$,
implies
no
natural 
Type-S 
generalized arithmetic $(\alpha,d)$
 can recognize  its
own  consistency
when $d$ is one of 
Example \ref{ex-2.1}'s three
  Hilbert-style 
%ss% versions of deduction:
deduction operators of 
$\, d_E \,$, $\, d_H \,$  
or
 $\, d_M ~~$:
\begin{quote}
{\bf ++ }
%\footnotesize
{\it 
(Solovay's  
modification
%Generalization 
\cite{So94}
%1994 Generalization \cite{So94}
%of a 1985 theorem 
of Pudl\'{a}k \cite{Pu85}'s formalism 
with
%using 
%some of 
Nelson and Wilkie-Paris \cite{Ne86,WP87}'s
methods)} :
Let 
$ \, (\alpha,d) \, $ 
denote 
a generalized arithmetic
supporting
% which contains 
the
\el{totdefxs}'s
Type-S statement
and 
assuring
the successor operation
will
satisfy
both 
% the axioms of 
 $  \,   x'     \neq 0     $ and
$     x'     =     y' \Leftrightarrow x=y $.
$~$Then
$ \, (\alpha,d) \, $  
%%%%$~\alpha~$
cannot verify its own
%will be unable to recognize its
%own  
consistency
whenever
simultaneously
 $d$ is
a Hilbert-style 
deductive
apparatus and
%whenever
$~\alpha~$
 treats addition and multiplication
as 3-way relations, 
satisfying 
their usual % identity,
associative, commutative 
 distributive 
and identity 
%  axiomatic properties.
axioms.
%   -axiom
% properties.
\end{quote}
Essentially, Solovay \cite{So94} 
privately communicated 
to us 
in 1994
%to us
an analog of theorem $++$.
%but 
Many authors
have noted Solovay
 has 
been
%often
reluctant to publish
% several of 
his 
nice 
privately communicated
results 
on many occasions
%in several contexts
\cite{BI95,HP91,Ne86,PD83,Pu85,WP87}. 
Thus,
%polished
approximate  analogs of 
%statement
 $++$
 were  explored 
subsequently
 by  Buss-Ignjatovic,
H\'{a}jek 
and
\v{S}vejdar in \cite{BI95,Ha7,Sv7},
as well as in Appendix A of 
our paper
\cite{ww1}.
Also, 
Pudl\'{a}k's initial 1985 article  \cite{Pu85} 
% implicitly
captured
% , notably, 
the majority 
%most
%%%  much 
of $++$'s 
essence, chronologically before Solovay's observations.
Also,
% and
Friedman did
% some
related work
 in
\cite{Fr79a}.

\item[   $~~$b.$~~$  ]
Part of what makes
%ss% the Pudl\'{a}k-Solovay discovery in
 $++$ interesting is that 
\cite{ww1,ww5,wwapal}
%Willard
developed two 
% separate
methods for 
% basis systems
generalized arithmetics
%%% $\alpha$ 
to confirm their own consistency, whose
natural hybridizations  are precluded by $++$.
Specifically,  these results involve
either a Type-NS
% basis 
system 
\cite{ww1,wwapal}
 verifying its own consistency
under 
any of our three
%main 
discussed
 variants of
% any  of the 
% $d_E$ or  $d_H$
% or $d_M$'s
Hilbert-style methods,
or   a Type-A 
%basis 
system \cite{ww93,ww1,ww5,ww6,ww14} 
verifying
 its
% own 
self-consistency
under $d_F$'s tableaux 
%deductive 
apparatus.
Also, Willard \cite{ww2,ww7} observed how one could
refine $++$ with Adamowicz-Zbierski's
methodology \cite{Ad2,AZ1} to show 
 Type-M  systems
cannot recognize their semantic tableaux consistency.
\ennd
\end{exx}

The roles of
% Observations
Items (a) and (b)
in our research
%from the current example, 
will become more evident
as this article progresses. 
Essentially, our
prior research
% ,
% best summarized in \cite{ww14},
%  has 
had focused mostly on Type-A arithmetics
that could verify their consistency under either semantic
tableaux deduction or some near-cousin of this concept
(as was explained in \cite{ww14}'s short 16-page summary
of \cite{ww93}-\cite{ww9}'s results).
The main challenge, faced by our prior research, was that $~ ++ \,$ 's
generalization of the Second Incompleteness Theorem for Type-S
arithmetics
% forced 
caused
us to rely 
%%  mostly
on analogs of \cite{wwapal}'s  ``ISCE''
Type-NS framework,
when a logic recognizes its
  own Hilbert consistency. (The latter, unfortunately,
dropped the assumption that successor is a total function. It instead
constructed the infinite set of integers by employing
an infinite set of distinct constant symbols $C_1,~C_2,~C_3,~....~$
where $C_j\, = \, 2^{j-1}~.$) % What we shall do in 
In
this article,
% is show
we will improve upon ISCE by showing
how one needs only three built-in
constant symbols, corresponding to the integers of
0, 1 and 2, to construct an
``IQFS''
 axiom system that can verify its own Hilbert 
consistency $~---~$ where here 
IQFS has access to a new type of `` $\theta$ '' function
primitive 
(defined later)
for constructing the infinite range of integers
from these three starting constants.
% of 0, 1 and 2. 

\begin{deff}
\label{def-2.4}
\rm
Let
$~\alpha~$ again 
denote an axiom basis, 
and $~d~$ 
designate
 a
deduction apparatus.
Then 
the  ordered pair
 $~(  \alpha  , d  )$
will
% shall
be called  a {\bf Self Justifying}
configuration
 when:
\begin{description}
% \xxitch
% \small
  \item[  i   ] one of  $~(  \alpha  , d  )$'s  theorems
(or possibly one of $\alpha$'s axioms)
do
%will
state that the deduction method $ \, d, \, $ applied to the
basis
system $ \, \alpha, \, $ 
%will 
produces a consistent set of theorems, and
\item[  ii   ]
     the axiom system $ \, \alpha  \, $ is in fact consistent.
\end{description}
\end{deff}

\begin{exx}
\label{ex-2.5}
\rm
Using 
Definition \ref{def-2.4}'s
 notation, our
prior
 research  
in
\cite{ww93,ww1,ww5,wwapal,ww9,ww14}
developed
% ordered pairs  
arithmetics
$~(  \alpha  , d  )$
that 
 were
%are
``Self Justifying''.
It 
%  has 
also 
% explored how 
proved
the Second Incompleteness Theorem 
implies specific
%formalizes
limits beyond which 
self-justifying
%such 
% particular
formalisms cannot transgress.
For any  $\,(\alpha,d) \,$, 
it is 
% easy
%actually 
almost trivial 
to construct a 
% second
%%% axiom 
system $ \, \alpha^d \, \supseteq  \,  \alpha  \, $
 that  satisfies
the
Part-i 
condition
(in an isolated context {\it where the Part-ii condition is
% also
 not also
satisfied}).
%of 
% this definition.
For instance,  $ \, \alpha^d \, $  could
consist of all of $~\alpha \,$'s axioms plus 
an added {\bf $\,$``SelfRef$(\alpha,d)$''$\,$} sentence,
defined as stating:
%%%
%%% the following 
%%% %%% added
%%% further
%%% sentence,  
%%% called
%%% %%% that we call
%%% {\bf SelfRef$(\alpha,d)~$}:
\begin{quote} 
%\xxitch
$\oplus~~~$ 
There is no proof 
(using 
$d$'s deduction method)
of  $0=1$
from the  {\it union}
 of
the
 axiom system $\, \alpha \, $
with {\it this}
sentence  ``SelfRef$(\alpha,d) \,$'' (looking at itself).
\end{quote}
Kleene 
discussed
 in
\cite{Kl38} 
how
to
% roughly
encode rough
% essential
% approximate
 analogs of 
the above
%this
% {\bf $\,$``SelfRef$\,$ Declaration''}.
 {\bf $\,$``I Am Consistent''} 
axiom
%axiomatic
statement.
% declaration.
%%% SelfRef$(\alpha,d)$'s
%%% self-referential statement.
Each of
Kleene, 
Rogers and Jeroslow 
 \cite{Kl38,Ro67,Je71},
however,
%{\it
 emphasized
that
$\alpha ^d$ 
may
be inconsistent
(e.g.  violating Part-ii of   self-justification's
definition),
{\it despite SelfRef$(\alpha,d)$'s 
formalized
assertion.}
%%%It thus often
%%% violates Part-ii of   self-justification's
%%%definition.
This is because if the 
%  ordered
 pair $(\alpha,d)$ is too strong
then a classic G\"{o}del-style diagonalization argument can
be applied to the axiom system 
$\alpha^d~~=~~ \alpha \, + \, $ SelfRef$(\alpha,d)$,
where the added presence of the statement 
SelfRef$(\alpha,d)$ 
will cause this extended version of 
$\, \alpha\,$, ironically,
 to
 become automatically inconsistent.
Thus, the
encoding of
% machinery of the sentence
``SelfRef$(\alpha,d)$'' is relatively easy,
% to encode,
%make well-defined 
via an application of the Fixed Point Theorem,
but this sentence
is
ironically 
%%%%%{\it most often  
{\it 
typically}
%  entirely
%usually
useless! 
\end{exx}

%\newpage

Unlike our earlier work, which focused 
 mostly on  a 
semantic
tableaux apparatus for deduction,
the current paper
will
explore
%paper will explore
\dfx{def-2.2}'s
more pristine Hilbert-style methodologies.
%% 
%% analogous to 
%% Example
%% \ref{ex-2.1}'s 
%% textbook
%% methods.
%% 
% of 
% $d_E$,  $d_M$
% and  $d_H$.
%%! 
%%! in
%%! the textbooks by 
%%!  Enderton,  H\'{a}jek-Pudl\'{a}k,
%%! Mendelson and Papadamiriyou \cite{End,HP91,Mend,Papa}.
There are, of course, many types of generalizations
of the Second Incompleteness Theorem known to
arise in  Hilbert-like  settings
\cite{Be95,Be97,Be3,BS76,Bu86,BI95,Fe60,Fr79a,Go31,Ha7,Ha11,HP91,HB39,Lo55,Kr87,Kr95,Pa71,Pa72,Pu85,Pu96,So94,Sv7,Vi5,WP87,ww1,wwapal}.
Each such
generalization
formalizes
a paradigm where 
self-justification is infeasible
under a Hilbert-style apparatus.

\smallskip

Our 
main
prior research about
%%! 
%%! main work about arithmetics displaying knowledge
%%! about their 
%%! 
Hilbert consistency appeared in \cite{wwapal}.
Its ISCE$(\beta)$ 
formalism could recognize its own 
Hilbert consistency and 
%%could 
prove analogs of any 
r.e. 
extension $~\beta$ of
Peano Arithmetic's $\Pi_1$ theorems.
It unfortunately required
the use of an infinite number of constant symbols, with
 ISCE$(\beta)$ 
using one built-in constant
symbol 
$~C_i~$,
for each power of 2.
% An alternative in  \cite{wwapal},
% called ISINF$(\beta)$, required use of only three constant symbols,
% but its proof lengths were impractically long. 

%%! required
%%! % in excess of
%%! an impractical
%%!  $O(N)$ length proof to construct an integer $N$.

\smallskip

Prior to \cite{wwapal}'s publication, 
Pavel Pudl\'{a}k
\cite{Pupriv}
examined our article and asked 
% us
%%% the question about 
%crucial 
whether
one  could improve 
upon ISCE's properties
 by using Ajtai's observations
\cite{Aj94}
 about
the Pigeon Hole principle.
In a
briefer and
 more informal respect,
% and abbreviated sense, 
Sam Buss
asked us a similar question during a lunch at a 1997
logic conference \cite{Bu97}.
%%  
%% Sam Buss \cite{Bupriv} also asked us
%% a 
%% %similar 
%% related
%% question
%% (during a more 
%% informal 
%% %abbreviated
%% conversation).
%% %(in a more informal manner).
%% 
Our prior
partial
answer to  
these
%Pudl\'{a}k's
questions
was offered  
%issue
%appeared
%in Sections 6 and 7 of \cite{wwapal}.
in Sections 6 of \cite{wwapal}.
A different
% type of 
and
% much 
substantially
more interesting
reply 
to these questions
will be offered
% 
% We will offer 
% % an alternate much  
% a 
% % much
% more sophisticated
% and different 
% type of reply
% % analysis
% %  formalism 
% 
in
the current 
paper.

\begin{deff}
\label{def-2.6}
\rm
Relative to any fixed given axiom system $\, \gamma \,$,
% a primitive
a formal
function symbol
 $~F~$ 
% shall 
will
be called a 
{\bf Q-Function} 
iff 
$\, \gamma \,$
is 
sufficiently ambiguous
for there to exist an UNCOUNTABLY infinite number of
different 
{\it 
plausible  sequences} of
formally 
enumerated
ordered pairs 
$\,(\, i \, , \, a_i\,)\,$
in expression \eq{wow} where
$~F(i)=a_i~$  is allowed as a
% logically 
permissible
%plausible 
%formalization 
representation 
% of
for 
$F$
under 
$ \gamma$'s axioms.
\begin{equation}
\label{wow}
  (0,a_0) 
 ~,~  (1,a_1)   ~,~  (2,a_2)   ~,~  (3,a_3)   ~,~  (4,a_4)~ ...
\end{equation} 
\end{deff}

\gvs

\tttc
 
Most
Q-Function symbols are
unsuitable for 
analyzing
%producing a positive resolution to
Hilbert's  Second Open Question or most 
issues in
% other prominent
% % mathematical 
% questions within 
mathematics.
 This is because the
presence of an 
 uncountably
 infinite
  number of
different 
plausible  sequences,
using Line
\eq{wow} to specify
$~F(i)=a_i~,~$  is 
%typically 
more of a burden than a benefit
(in most cases).
An exception to this general rule
of thumb
 will be
provided by
 the next
section's $~\theta~$ operator.
It will be germane to 
 $\, ++ \,$'s generalization of the Second Incompleteness
and suggest a mechanism whereby an efficient form of
``Type-NS''self-justifying 
arithmetic
can recognize its own Hilbert consistency
in a 
%nontrivial 
quite substantial
sense.

% \section{Revisiting a World which Hilbert called
% {\it ``Cantor's Paradise''}}

\section{Starting Notation Conventions}
\label{ss4}
\label{seee3}

%333333333333333333333333333
\vspace*{- 0.4 em} 

% OLD Title was {\it Notation and Basic Concepts}

\tttc

% Throughout this paper,
%article,
A function 
 $\, H \, $ 
will be     called
 {\bf Non-Growth} 
whenever
%iff  
$ H(a_1,a_2,...,a_j) 
\leq   Maximum(a_1,a_2,...,a_j)$
for all  $a_1,a_2,...,a_j$.  Six  examples of  
 non-growth functions are:
\bee
\small
% \baselineskip = 1.05 \normalbaselineskip 
\topsep -5pt
\itemsep -1pt
\parskip 0pt
\ttt2c
\item
{\it Integer Subtraction} 
where ``$~x-y~$'' is defined to equal zero 
%ss% in {the special case} where
when $~x \leq y,$
\item
{\it Integer 
Division}
where ``$~x \div y~$'' equals
$~0~$ when $~y=0~$ and
it equals $~\lfloor ~x/y ~\rfloor~$ otherwise,
\item
$Maximum(x,y),~~$
\item
$Root(x,y)~$ which equals $~ \lfloor ~x^{1/y}~ \rfloor$ when $~y\geq 1~$
%% 
%% and
%% it equals $~x~$ when $~y=0.$
%% 
(and zero otherwise),
\item
$ Logarithm(x)~=~\lfloor ~$Log$_2(x)~ \rfloor~$
and
\item
$Count(x,j)~~=~~$the number of ``1'' bits
among $~x$'s rightmost $~j~$ bits.
\ene
%% 
%% 
%% \bee
%% \baselineskip = 0.8 \normalbaselineskip 
%% 
%% \item
%% {\it Integer Subtraction} 
%% where ``$~x-y~$'' is defined to equal zero 
%% in {\it the special case} where
%%  $~x \leq y,$
%% \item
%% {\it Integer 
%% Division}
%% where ``$~x \div y~$'' equals
%% $~x~$ when $~y=0,~$ and
%% it equals $~\lfloor ~x/y ~\rfloor~$ otherwise,
%% \item
%% $Root(x,y)~$ which equals $~ \ulcorner ~x^{1/y}~ \urcorner$ when $~y\geq 1,~$
%% and
%% it equals $~x~$ when $~y=0.$
%% \item
%% $Maximum(x,y),~~$
%% \item
%% $ Logarithm(x)~=~\lfloor ~$Log$_2(x)~ \rfloor~$ when $~x \geq 2,~$
%% and  zero otherwise. 
%% \item
%% $Count(x,j)~~=~~$the number of ``1'' bits
%% among $~x$'s rightmost $~j~$ bits.
%% \end{enumerate}
%% 
These operations were called
either the 
{\bf ``Grounding''} 
or {\bf ``Ground-Level''} 
 functions
in our articles \cite{ww1,wwlogos,ww5,ww14}.
We will
use the latter nomenclature in the current article
 because the notion of a  ``Ground-Level''
function should not be confused with the {\it  very different} notion
of a ``Grounded Term'' employed by Definition
\ref{def-3.4}. 

%  TWO DEFS or ONE ???????? 

Our
starting  language $L^G$  
%shall
will
 also contain
the
two atomic 
symbols
%% relations 
of ``$~=~$'' and ``$~\leq~$'' and three
built in constants symbols, $~C_0~$, $~C_1~$ and $~C_2~$, 
for representing
the values of 0, 1 and 2. 
%Within 
In
this context, expressions 
\eq{newadd} and \eq{newmult} formalize 
exactly
how addition and multiplication
can be encoded as two 3-way predicates,
%%  in  $L^G$,
 denoted as
Add$(x,y,z)$ and Mult$(x,y,z)$.

\vspace*{- 0.6 em} 
\beq 
\label{newadd}
z ~ -~x~~=~~ y~~~~ \wedge ~~~~ z~\geq~x
\end{equation}

\vspace*{- 1.2 em} 
\begin{equation}
\label{newmult}
[~(x=0    \vee    y=0 ) \Rightarrow z=0~ ]~ ~\wedge ~~ 
[~(x \neq 0 \wedge y \neq 0~) ~ \Rightarrow ~
(~ \frac{z}{x}=y  ~\wedge \, ~  \frac{z-1}{x}<y~~)~]
\end{equation}
Also, 
our constant symbols,
 $~C_1~$ and $~C_2~,~$ for representing
the quantities 1 and 2, 
 allow us to
formalize the following further
non-growth
%useful 
%%%%%%%% function 
operations:
\bed
\small
\topsep -5pt
 \itemsep -2pt
\ttt2c
\it
\item[   $~$a.  ] 
Pred$(x)~~=~~x-1$
(in a context where 
%Item 1's
the prior paragraph's
%our 
definition for Subtraction implies
Pred(0)$=0$.$~)$   
\item[   $~$b.  ] 
Half$(x)~~=~~x \div 2$
(in a context 
where 
``$ ~x \div 2~$'' equals technically 
$~\lfloor \, \frac{x}{2} \, \rfloor~$ 
under 
%our
the 
%preceding 
prior
paragraph's 
notation).
% of Item 2's definition for Division).  
\item[   $~$c.  ] 
Pred$^n(x)$ defined to be $~n~$ iterations of the
Predecessor operation
\item[   $~$d.  ] 
Half$^{ \,n}(x)$ defined to be $~n~$ iterations of the
halving operation.
%%
%%\item[   $~$e.  ] 
%%Bit$(x,j)~=~$ $Count(x,j)~-~Count(x,j-1)~$
%%(e.g. the value of $~x \,$'s rightmost
%% $j-$th bit)
%%\item[   $~$f.  ] 
%%Min$(x,y)~~=~~x~-~(y-x)~~$
%%(e.g. a quantity that represents the minimum of $x$ and $y$
%%because our definition of ``integer subtraction'' 
%%implies  $~y-x \,= \, 0 ~$ whenever $\, x \geq y$.$~~)$  
%%
\ennd

Let us say
that
 a function symbol $H(x_1,x_2...x_j)$
is
 {\bf Growth Permitting}
iff for each integer $~k \geq 2~$ there exists a 
``growth-tuple''
$(a_1,a_2,..,a_j)$ 
satisfying
  $~H(a_1,a_2,..,a_j) \, >\, k~$
and 
also
% simultaneously 
having
each $ \, a_i \leq k$ \, 
It will  be necessary 
%% for  us 
to employ
 either an infinite number
of constant symbols or some Growth-Permitting function, 
so that an extension of
our
% the 
language $L^G$ can construct the 
%infinite
%full collection of
naturally
extended
%infinite
broader set
% broader range 
of 
integers 
%of 
$~3,4,5,6~....~$ from its starting objects of 0,1 and 2.

\smallskip

One method for resolving this problem has
been presented in \cite{wwapal}.
Its  ISCE$(\beta)$ axiom basis did
 employ  an infinite number of further constant symbols. It
was
compatible with                     self-justification,
but deviated from 
%{\it very sharply from} 
Hilbert's
intended
 goals 
because it employed
% by employing 
%an  
a {\it highly awkward}
infinite number of 
distinct
 constant symbols.

\smallskip

The challenge posed by $++$ 
is, thus, formidable.
Our goal in 
%the current 
this
article
will be to suggest 
how
%that 
a Q-function primitive $F(x)$,
constrained by
deliberately
ambiguous
axioms,
% function definition,
can help overcome the constraint that $++$ imposes. 
Such an
% ultra-ambiguous defined 
unusual
primitive $~F~$ will have
an  uncountable number
of vectors, analogous to Line \eq{wow},
%  that satisfy
as permissible solutions for
$F$'s
definition.
Our basic goal
% in this article 
will be to outline 
how this unusual concept is
% likely  
germane to the
self-justifying
%% 
%% why it is
%% likely such Q-functions will enable one to build
%% surprisingly efficient self-justifying logics that
%% partially (but not fully) achieve the 
%% %{\it  diluted portion} of the
%% 
aspirations  
 that
Hilbert and G\"{o}del
had
expressed in 
%statements 
$*$ and
$**~$.

% \smallskip

\begin{deff}
\label{def-3.1}
\rm
Let us say 
% an integer-valued
a 
formal
function symbol
$F(x_1,x_2,...x_j)$ is {\bf ``1-Definitive''} iff it has only one
permissible
solution under 
an axiom system $\gamma$'s
% formal 
definition 
of it.
Let us call
the function $~F~$  {\bf ``Indeterminate''} otherwise.
%(The remainder of this article 
\end{deff}

% \medskip

Mathematicians
obviously typically avoid using 
axiom systems that entail using ambiguous function
definitions (e.g. they certainly usually avoid analogs of
\el{wow}'s 
dizzying
quantity
of possibly
$\aleph_1$ distinct  
solutions for its defined function
% symbol
 $F$ ).
% 
% function definitions
% %%%%%%%%%, $~F~$ 
% with 
% %that have 
% %%%%  having
% even two solutions, not to
% speak of
% what will be
% \el{wow}'s 
% dizzying
% quantity
% of possibly
% %potential
% %possible
% $\aleph_1$ distinct 
% % number of 
% solutions,
%
%considered in the current article.
%
Our effort to overcome 
$++$'s 
broad-scale
generalization of the 
Second 
Incompleteness
Theorem will be unconventional because we will
use an indeterminate function,
%  symbol, 
called the
$~\theta~$ operator, 
to overcome $++$'s 
challenge. It will turn out that the 
$~\theta~$ 
operator will almost duplicate the logarithmic
efficiency of the traditional addition and multiplication function
symbols for constructing the natural numbers, while simultaneously
dove-tail around $++$'s generalization of the Second Incompleteness
Theorem. 

The  indeterminate nature of 
$\, \theta\,$ will be crucial to our evasion of the Second
Incompleteness Effect.
Our research into this topic was
greatly
 influenced
by an email we received from 
Pavel
Pudl\'{a}k
\cite{Pupriv}, 
shortly
after he received an advanced copy of our
article 
 \cite{wwapal}.
He
appreciated the nature of the challenge 
the
% \cite{wwapal}'s
 ISCE axiom framework faced,
when it used  an awkwardly infinite number of constant symbols.
%% Pudl\'{a}k appreciated
 (ISCE relied on
% using 
this
% unfortunate 
infinite magnitude to evade
%invoking
 $++$'s 
requirement
%invariant
% result 
that no evasion of the
Second Incompleteness Theorem
can take place
 when an
arithmetic recognizes 
 Successor as a total
function.)
Pudl\'{a}k's  
%private 
%His
emailed communications
\cite{Pupriv}, thereby
suggested 
% that 
we look at 
Ajtai's
work 
\cite{Aj94}
about a
%the 
Pigeon-Hole function
 $~ \glamb(x)\,$, defined by the identities
\eq{zm1} and \eq{zm2}.

%\newpage
% \vspace*{- 1.2 em} 
\beq
%% \small
\label{zm1}
\forall ~~x~~~~~ \glamb(x)~ \neq ~ 0
\enq

\vspace*{- 1.2 em} 

\beq
\label{zm2}
%% \small
\forall ~~x~~~ \forall ~~y~~~~ x ~ \neq~ y ~~ \Rightarrow ~~
\glamb(x)~ \neq ~\glamb(y)
\enq
The relevance of 
$~\glamb~$ 
% Pigeon-Hole functions
can be
best 
%readily 
appreciated
% if 
when
%we let
$~\glamb^n(x)~$ 
 denotes
% the 
a
term
 $~\glamb(~\glamb(~ ... \glamb(x)))~$
consisting of $~n~$ iterations of the $~\glamb~$ operator.
% Then 
% the
\el{DUMB1}'s
composite
term $~S_n~$
% , defined below, shall
will then satisfy  
the lower bound of $~S_n \,  > \,n$. 

%% An axiom system, employing the primitive 
%% operation
%%  $~ \glamb~,~$ 
%% can thus 
%% can easily
%% prove
%% Line \eq{farreach}'s 
%% assertion.
%% %claim.
%% %under almost all conventional logics.
%% 
%% 
\beq
\label{DUMB1}
S_n~~~=~~~\mbox{Max}[~\glamb(0)~,~\glamb^2(0)~,~\glamb^3(0)~,~...~~\glamb^n(0)~]
\enq
Pudl\'{a}k
observed
% that
%the
% Pigeon-Hole function 
 $~ \glamb(x)~$  
will
grow too slowly (in the worst case)
for
one to be able to 
deduce
successor is a total function
from its properties  
% further observed that it is known 
 \footnote{ \baselineskip = 0.94 \normalbaselineskip
 The operation $\glamb(x)$ will grow
at a slower rate than Successor, 
if it equals $x+1$ for all standard
numbers $~x~$ and if $\glamb(x)=x-1$ 
 when $~x~$ is
a non-standard integer. This seemingly minute detail
implies  one cannot infer 
Successor is a total function from
 $\glamb$'s behavior, since the latter is contradicted by a
 model  where
 all  non-standard
numbers have 
%their 
sizes bounded by some fixed  
% non-standard 
number B.
(This 
subtle detail, 
raised by
Pudl\'{a}k's email \cite{Pupriv}, was fascinating because
it 
%shows that
raised the question about whether
 a partial exception to
Example \ref{ex-2.3}-a's
invariant $++$
%% on \pag2,  
might plausibly exist.)  }.
%% 
%% thus, 
%% suggests the 
%% Pudl\'{a}k-Solovay 
%% version of the Second Incompleteness
%% Theorem (stated on \pag2) 
%% might
%% %%%%%should
%% allow for
%% potential
%%  exceptions 
%% to it
%% arising from the 
%% %delicate 
%% formal
%% behaviour of
%% some
%% %% 
%% %% presence of
%% %% %some
%% %% these permissible
%% %% 
%% %% 
%%  non-standard 
%% variants of
%% % interpretations for 
%% the  Pigeon-Hole function $\glamb$.  }.
%% 
%% 
%that 
%prove 
His 
%insightful 
email \cite{Pupriv}
% thus
 hinted
% that
% asked whether 
the inequality 
$S_n \, >  \, n \,$
might
%would,
% thus, 
% still
enable a formalism,
%,
% based around
utilizing the
 $\, \glamb \,$ operative, 
to 
% somehow
improve upon \cite{wwapal}'s results.

% our
% formalisms could be 
% revised
% %modified 
% so that 
% % the  Pigeon-Hole function 
%  $~ \glamb(x)~$ 
%  could improve upon \cite{wwapal}'s results.

%%
%%(possibly using Ajtai's methodologies \cite{Aj-focs}).
%%Sam Buss raised, interestingly,  a 
%%partially
%%similar 
%%issue during an informal conversation
%% \cite{Bu-priv} in 1977.
%%
%%\smallskip
%%
%%These questions
%%% by
%%%Pudl\'{a}k and Buss 
%%were insightful because they isolated 
%%an
%%important juncture where $++$'s underlying methodology does not apply.
%%A partial answer to these questions appeared in 
%%\cite{wwapal}'s closing section, but a more comprehensive full
%%answer  has always eluded us.

%This is  because there always seemed to appear
%one wrinkle of details that precluded a full proof.

\smallskip

It  was
initially
 unclear 
%%%%% to us
whether a positive answer to 
Pudl\'{a}k's
interesting
% probing
 question would resolve ISCE's main difficulties.
This is
because
% Expression 
\eq{DUMB1}'s
term
$~S_n~$  requires $O(~n^2~)$ logic symbols to encode
% essentially 
an integer quantity
greater than
 $~n~$
(since its term
$~\glamb^j(0)~$ uses $O(j)$ logic symbols).
%an integer quantity that exceeds the quantity $~n~$ in size.
Thus, the integer $~2^{100}~,~$ whose binary encoding
%used 
uses
100 bits,  would
actually
 require 
% fully
in excess of 
 $~2^{200}~$ bits to encode. 
Such quantities, exceeding the number of atoms in the universe,
were troubling because our 
general
goal has been to 
construct self-justifying arithmetics that 
 possessed, at least,
some
partial  facets of
 pragmatic value.

% 
% find a partial
% answer to Hilbert's
% Year-1900 Second
% Problem  
%  that would 
%  possess, at least,
% some
% partial  facets of
%  pragmatic value.
% 

\bigskip
\gvs
\tttc

The remainder of this section will outline how a different type of
Q-Function operator will 
be much better than
 $~ \glamb~$ for meeting our needs.  
During our discussion,
Power$(x)$ will denote 
a primitive specifying
% that
 $~x~$ is
a power of
$~2~$.  It is 
%formally 
encoded
by 
\eq{wep2} 
because
%under
our Grounding language
has
``Logarithm$(x) \,$''$ ~ = ~ \lfloor \,$Log$_2(x) \, \rfloor \,$.
\beq
\vspace*{- 0.6 em}
\label{wep2} 
%\small
x=1 ~~~\vee ~~~ \mbox{Logarithm}(~x~)~\neq~\mbox{Logarithm}(~x-1~)
\enq
In this context,
 $\zzthe(x)$ 
will denote 
%the
our new analog of
the $\glamb(x)$ function  
%% haphazard
that walks among the powers of 2 in a manner 
similar to 
$\glamb(x)$'s
%  haphazard
 walk through conventional
integers. 
It is 
% formally
defined by \eq{walk1}-\eq{walk4}.
% 
% It will thus satisfy
% the axiomatic constraints below (which are 
% $\zzthe(x)$'s analog of the more modest constraints given in
% % sentences 
% \eq{zm1} and \eq{zm2}).
% The most important difference between these two constructs
% is that axiom \eq{walk1} requires that 
%  $\zzthe(x)$ maps power of 2 onto powers of 2.

% {\small

\vspace*{- 0.6 em}
{\parskip -6 pt
\beq
\label{walk1}
\forall ~~x~~~~~ \mbox{Power}(x) ~~~ \Rightarrow  ~~
\mbox{Power}(~ \zzthe(x)~)
\enq
\beq
\label{walk2}
\forall ~~x~~~~~ \zzthe(x)~ \neq ~ 1
\enq
\beq
\label{walk3}
\forall ~~x~~~ \forall ~~y ~~~~[~ x ~ \neq~ y ~
\wedge ~\mbox{Power}(x)~]
 \Rightarrow ~~ 
\zzthe(x)~ \neq ~\zzthe(y)
\enq
\beq
\label{walk4}
\forall ~~x~~~~~ \neg ~ ~\mbox{Power}(x)~~~~ \Rightarrow ~~~~ 
\zzthe(x)~=~0
\enq}

\vspace*{- 1.2 em}
\noindent
{\it It needs to be emphasized} that 
 \eq{walk1} -- \eq{walk4} will be the
{\it only vehicle} our
self-justifying axioms 
%will 
have available to construct
integers $\geq \, 3$. $~$These
axioms will be called
%They will be henceforth called
the {\bf $~$Up-Walking$~$} axioms. 
(The axiom \eq{walk4} 
is,
% does not,
technically, 
unnecessary
 to construct 
% any
integers  $\geq \, 1\, $, but it is helpful because
it 
% allows us to 
formalizes how our methodology will treat  integers
which are not powers of 2.)

Both
the 
Q-functions 
% the operators 
$\glamb$ 
 and $\zzthe$ are 
challenging
%daunting
because there are a  
 dizzying
$\aleph_1$ distinct 
vectors, analogous to 
 Line \eq{wow}, 
that are
%where their definitions permits 
representations of these functions.
%% 
%% Also, we may combine either operation with our 
%% language $L^G$'s grounding function-primitives to formulate a term
%% $~T_n~$ that defines any arbitrary integer $~n~$.
%% 
We will soon see,
 however,  that
there is
a
sharp
distinction
%  major difference 
between these
% two 
concepts
from a computational
 complexity 
perspective.

\begin{definition}
\label{defx-3.2}
\rm
Let $~L^Q~$ 
and $~L^{Q^*}~$ 
denote the 
extensions
of $~L^G\,$'s Grounding language that contain the
respective
additional 
function symbols of  
  $\zzthe$
 and
$\glamb$. Then 
$~~L^Q~$ shall be called the {\bf  Q-Grounding} language, and 
 $~~L^{Q^*}~$ 
will be called the  {\bf  Q* Grounding} language. 
\end{definition}

\begin{propp}
\label{th-3.3}
In contrast to the
Q* Grounding language 
that requires $O(~n^2~)$ function symbols
for defining a term $~T^*_n~$ for representing the integer
$~n,~$ the    Q-Grounding language
%% will need no more than 
needs
% uses 
only
$O(~$Log$^{ \, 3\,} \,n~)$ symbols to 
encode
%formalize 
a term 
$~T_n~$ representing
$~n$.
\end{propp}

\vspace*{- 1.0 em}

\begin{center}
% \small
% Our proof of \phx{th-3.3} 
\phx{th-3.3}'s 
proof
will rely upon the following notation  convention:
\end{center}

\vspace*{- 0.8 em}

\begin{definition}
\label{def-3.3}
\rm
Let
 $~\zzthe^j(x)~$
denote the term
 $~\zzthe(~\zzthe(~ ... \zzthe(x)))~$
where there are 
$~j~$ iterations of the 
 $~\zzthe~$ operation.
% Throughout this article,
%ss% Then
%for any  $~j \, \geq 1~,~$ 
%the symbol 
Let $~E_j~$ 
% shall
% will 
denote
the quantity produced by
\eq{ej-def}'s division operation: 

\vspace*{- 0.6 em}

\beq
\small
\label{ej-def}
 \frac{~\mbox{Max}
~[~\zzthe^{\, j \,}(1)~,~~\zzthe^{\, j-1 \, }(1)~,~... ~\zzthe(1)~~] }
{~~\mbox{Half}^{\,j\,} ~ \{ ~\mbox{Max}~[
 ~\zzthe^{\, j \,}(1)~,~~\zzthe^{\, j-1 \, }(1)~,~... ~\zzthe(1)~]~ 
 ~ \}~~ } 
%
% \mbox{Max}(~\zzthe^j(1),~\zzthe^{j-1}(1),~... ~\zzthe^1(1)~ 
\enq
It is easy to 
establish that
%see
$  E_j  =  2^j  $ for every
$j   \geq 1$.
This is
 because \eq{ej-def}'s
twice-repeating term
%ss% object
 of
$  \mbox{Max}
~[~\zzthe^j(1),  \zzthe^{j-1}(1),...\zzthe(1) \,]$
% is at least as large as $\, 2^j\,$.
is a power of 2 exceeding  $\, 2^j\,$.
%% 
%% The definitions of the
%% % Q-Grounding 
%% functions of ``Half'', ``Max'' and
%% ``$~\zzthe~$''   imply 
%% $~E_j~=~2^j~$ for each
%% $j \, \geq 1$. 
%% 
For the additional case where $~j=0~,~$
we 
shall 
%will simply
define  $~E_0~=~1~$ (by 
having the 
%%
%%setting it equal to 
%%our
%%%the
%%
formal constant symbol of
 $~C_1~$ stand for ``1'' ).
\end{definition}

%% , which 
%% is intended to
%% %formally 
%% represent the integer of ``1'').

{\bf Proof of  \phx{th-3.3}:}
%The justification of   \phx{th-3.3} is an 
Easy consequence of
\dfx{def-3.3}'s machinery. Thus if $~n~$ is a power of
2 of the form $~2^j~$ then
% the preceding
% definition's
expression $~E_j~$ is a term representing $~n \,$'s value
that employs
  $O(~$Log$^{ \, 2\,} \,n~)$ 
logical
symbols. On the other hand, if  
 $~n~$ is not a power of
2 then it can be defined 
with   $O(~$Log$^{ \, 3\,} \,n~)$ symbols by 
setting
$~E_j~$ equal to the least power of 2 greater than $~n~$ and
subtracting from $~E_j~$ those powers of 2 that are needed to
produce $\,n\,$'s  value. 
For example since $76~=~128~-~32~-~16~-~4~,\,$ it can
be formalized as a term $T_{76}$ defined by
$~E_7 \, - \,E_5 \, - \,E_4 \, - \,E_2 ~$.

% $~~~~\Box$

% \baselineskip = 1.8 \normalbaselineskip 

\begin{definition}
\label{def-3.4}
\rm
A ``term'' in mathematical logic
is defined to be a syntactic object,  built
out of
solely
 symbols for representing 
functions,
constants and variables.
% 
% The nomenclature in
% % classical
% logic has 
% %formally 
% defined
% a {\it ``term''} to be a syntactic object,  built
% out of symbols for representing 
% functions,
% constants and variables.
% 
% 
Such an object is called
% either 
a {\bf ``Ground Term''} 
%% (or for precision a
%%% {\bf ``Tree-Oriented Ground Term''} )
when it  is built {\it out of solely}
function and
 constant symbols.
For example in our Q-Grounding language (which 
uses
%owns only
$  C_0  $, $  C_1  $ and $  C_2  $ 
as
built-in
 constants),
% symbols), 
the expression
%of 
``$\, C_2-  C_1\,$'' 
is
% such 
a ground term.
Two more complex 
examples of 
ground terms are
``Max$(   C_2 ,  C_1 - C_0)$'' 
and ``Max$( ~\zzthe(C_1)~,~C_2 ~)$''.
Also,
expression $~E_j~$ 
in Line \eq{ej-def}
should be viewed as 
a ground term (when one 
views
its
use of the
symbol
 ``1'' as an %informal 
abbreviation
for the
constant
 ``$~C_1~$''). 
\end{definition}

\begin{remm}
\label{rem-def-3.4}.
\rm
We will distinguish in
Proposition \ref{th-7.1} and
in
 other parts of
\textsection \ref{ss7} 
 between two
kinds of ground terms, that are called the
%
%{\bf Comment of Definition \ref{def-3.4}'s Notation:} 
%We will  distinguish  between two
%kinds of Ground terms in Section \textsection \ref{ss7},
%called its 
{\bf ``Tree-Oriented''} and 
{\bf ``Dag-Oriented''} formats.
The latter will differ from a more 
conventional tree structure
by having a 
Directed Acyclic Graph structure replace 
a 
%ss% logic's
term's
usual 
 tree format for defining its quantitative values.
Our discussion in the next several sections will be
simplified if we use the shorter phrase of
{\it ``Ground Term''} 
to refer to what
% Section 
\textsection \ref{ss7} will more
accurately called a 
{\it ``Tree-Oriented Ground Term''}.
(It will turn out
% that our   
Proposition \ref{th-7.1} 
will later explain how Dag-oriented ground terms
differ from
their
 tree-oriented 
counterparts
% Ground Terms by allowing us
%to reduce the 
by reducing the
$O(~$Log$^{ \, 3\,} \,n~)$ length of a
tree-oriented term to a more compact
$O(~$Log$\, \,n~)$ size.)
 \end{remm}

\begin{definition}
\label{def-3.5}
\rm
A ground term $~T~$ will be called an 
{\bf ``Observable''} 
object iff there is
%{\it only one} 
an
unique
interpretation of its 
quantified value in the
%meaning  in our
Q-Grounding language.
It 
%will be
 is
called an
{\bf ``Unobservable''} iff it has multiple 
%plausible 
such
interpretations
due to $\zzthe$'s ``indeterminate'' definition
(e.g. see Definition \ref{def-3.1}). 
\end{definition}

%%% (due to the 
%%% %uncountably 
%%% ambiguous nature of 
%%% % our built-in function 
%%% $~\zzthe~~$).
%%% \end{definition}

\begin{exx}
\label{ex-3.6}
\rm
The previously mentioned ground term
Max$( ~\zzthe(C_1)~,~C_2 ~)$ is an ``unobservable''
because it can assume any of the plausible integer values
of $~2 \, , \, 4 \, , \,8  \,  , \,16  \, 
 \, ... ~$.
On the other hand, 
%% expression 
\eq{xoo} 
is an ``observable''
that
 represents
 the integer value of ``3''.
(This is because 
its 
twice-repeating
term 
``$~\mbox{Max}[ ~C_2~,  ~\zzthe(C_2)~, ~\zzthe^{ \, 2 \,}(C_2)~]~$'' is bounded
below by 4, causing the left and right sides of its subtraction
operation to differ by 
% an amount of 
exactly 3.) 
\beq
%% \small
\label{xoo}
\mbox{Max}[ ~C_2~, ~\zzthe(C_2)~, ~\zzthe^{ \, 2 \,}(C_2)~] ~~-~~
\mbox{Pred}^{\, 3 \,} \{~\mbox{Max}[ ~C_2~, ~\zzthe(C_2)~, ~\zzthe^{ \, 2 \,}(C_2)~]~\} 
\enq
Our notation 
%thus 
also
implies that Line \eq{ej-def}'s
 expression $~E_j~$ 
is an
 ``observable''. This implies, in turn, that 
 \phx{th-3.3}'s term $~T_n~$ is an  ``observable''
 employing no more than 
 $O(~$Log$^{ \, 3\,} \,n~)$ 
logical
symbols.
For example since $76~=~128~-~32~-~16~-~4~,\,$ 
it follows that $~ T_{76}~$
corresponds to the
observable
 term
$~E_7 \, - \,E_5 \, - \,E_4 \, - \,E_2 ~$,
$~$where each $~E_j~$ employs only 
 $O(~$Log$^{ \, 2\,} \,j~)$ symbols.
%%%%%\end{exx}
\end{exx}

Thus, \dfx{def-3.5} and Example \ref{ex-3.6} have illustrated
%that 
how
the realm of ``observable'' objects is a 
% very 
broad and accessible world,
of 
non-trivial
%% pragmatic
 significance.
It allows every integer $~n~$ to be represented by a
% reasonably small 
term $~T_n~$ with
% an
a tight
 $O(~$Log$^{ \, 3\,} \,n~)$ length
(in a context where
Proposition \ref{th-7.1}'s
%  Section \textsection \ref{ss7}'s 
more elaborate formalism
will allow us to reduce this length to a 
yet 
more
attractive  $O(~$Log$ \,n~)$ size).

\medskip

The distinction between
``Observables'' and ``Unobservables''
% ground terms 
will 
%also 
%% 
%% cast a
%% delightful
%% 
offer a
 new perspective on
the aspirations
% that 
which
Hilbert and G\"{o}del 
stated
%expressed
in 
their
statements $*$ and $**$,
under 
%our
the 
new
%2-part conjecture
%(mentioned later in this paper).
IQFS axiomatic framework
that 
% will 
is
% shall be 
proposed later in this article.
It
 will suggest
how the Second Incompleteness Theorem
can
%  remain to 
be seen as a
fascinating
 {\it majestic result}
from a purist perspective, while
a {\it well-defined fragment} of
%their 
what
Hilbert and G\"{o}del
sought in 
 $*$ and $**$
%aspirations
%% in  statements $*$ and $**$
can 
% likely
%almost certainly
be 
%part-way
satisfied, in at least a 
% well-defined 
diluted
%limited
 sense.

\begin{remm}
\label{rem-3.7}
{\bf (explaining the goals of this paper):$~$}
\rm
Let us say 
%  that
a basis axiom system $~\alpha~$ owns
a {\it ``Finitized Perspective''}  of the Natural Numbers
if it requires only a 
{\it finite number} of proper axioms
to construct the full set of integers
$~0,1,2,3~... ~$. All conventional arithmetics have this property.
%It is useful to divide such 
Such
logics
%arithmetics
fall into two categories,
called {\it Single} and {\it Double-Formatted} systems.
%as defined below: 
They are defined below: 
 %These constructs are defined below:
\bed
\item[   a.  ]
%An axiom basis $~\alpha~$
%will be called
{\bf Single-Formatted Arithmetics} consist of
axiomatic basis systems
% $~\alpha~$
%all of 
whose 
%iff all its 
ground terms are 
all 
Observables.
(Most conventional arithmetics
%%%% will 
%fall into this 
lie in this
category
%when the 
because they typically employ
%%% employ the
%%%  growth
%%% %  function
%%% properties 
%%%  of
the
Successor
operation
%function
 in
% a straightforward 
%the
%% a  conventional
% the 
a
traditional
manner.)
%% since the 
%% the simple  growth function of Successor
%% easily
%% generates
%% all the natural numbers). 
%are {\it ``Single-Formatted Formatted''} logics. 
\item[   b.  ]
{\bf Double-Formatted Arithmetics} 
% representing 
represent
systems
%%%consisting of
%%%%axiomatic 
%%%logics
%%%%%basis systems
whose ground terms 
may be either
 Observables
or Unobservables.
(Axiomizations
for Q-Grounded logics 
%%  of
%% the
%% % our
%% Q-Grounding language
are
%%% will 
%obviously 
%%% be
``Double-Formatted''
because they
allow $\theta$'s analog of 
\el{wow}'s function symbol $F$ 
to have
an uncountable number of 
different allowed
representations).
% 
% (Our
% Q-Grounding language 
% gives support to such a system.
% This is because it
% can  have its function primitives
% defined by a finite number of 
% proper axioms,)
% %axiom-sentences.) 
\ennd
The distinction between 
categories
%Items 
(a) and (b) is
significant
% important
because 
Example \ref{ex-2.3}-a
%%%  \pag2 
%%% had 
%already 
explained how
 statement $++$'s generalization
of the Second Incompleteness Theorem applies to 
{\it any formalism} recognizing Successor as a total function.
Thus,
% Item (b)'s
 Double-Format logics 
%  Double-Formatted logics 
are useful, if one wishes to consider alternatives
where 
%formalism that do not recognize
successor
is not seen
 as a total function.
%More precisely, 
In this context,
Hilbert's 
%  famous 
%Year-1900 
Second 
Open
Problem
can be viewed
 as a  {\it  2-part question}, 
composed of sub-queries Q-1 and Q-2:
%%%%%
%%%%% {\it  2-part question}.
%%%%%The separation of Hilbert's question into two parts,
%%%%%called Q-1 and Q-2, will allow 
%%%%%%% 
%%%%%%% This 
%%%%%%% bipartite
%%%%%%% distinction 
%%%%%%% is useful because it
%%%%%%% can enable 
%%%%%%% 
%%%%%the academic community to better
%%%%%  with
%%%%% what Hilbert and G\"{o}del were
%%%%%seeking to accomplish
%%%%%in 
%%%%%their
%%%%%statements
%%%%%of $*$, $**$ and $***$.
\bed
\ttt2c
%%$\small
\item[ {\bf Question Q-1$~~$}] {\it Are any axiom systems
able to
 prove 
theorems 
verifying
 their own consistency in a robust sense?$~~$}
The answer to Q-1 is clearly ``No'' because the combination
 G\"{o}del's initial 1931 result \cite{Go31} with 
%the 
%further
Hilbert-Bernays's result 
\cite{HB39} 
and the Pudl\'{a}k-Solovay invariant $++$
(from Example \ref{ex-2.3}-a)
%% \pag2) 
imply 
arithmetics of ordinary strength cannot prove
their own consistency in a robust sense.
\item[ {\bf Question Q-2$~~$}]
 {\it Can 
logic systems
%arithmetic logics
%axiomatizations of Arithmetic
%  , at least, 
%somehow 
``appreciate'' 
% (not formally ``prove'')
 their
own consistency in some 
{\bf REDUCED} sense, that is diluted
but not
 fully 
immaterial?}
$~~\,$The answer to 
%question 
Q-2 is 
complex
%%% more complex than Q-1
%less clear-cut
because 
%several types of 
some
logics,
such as
our proposal
in \textsection \ref{ss5}
 and
 \cite{ww93,ww1,ww5,wwapal,ww9,ww14}'s weaker
and earlier paradigms,
 can
formalize
% ``recognize''
their 
own consistency, 
using Example \ref{ex-2.5}'s
% a 
Fixed-Point {\it ``I am consistent''}
axiom.
\ennd

A theme of this article will be that
% distinction 
the distinguishing
between questions Q-1 and Q-2 and between
Single and Double-Formatted Logics 
is 
related
% likely central 
to the mystery
% that has enshrouded 
enshrouding
the Second Incompleteness Theorem.
This is
%is germane to the aspirations of automated theorem proving
%will be germane to this article 
because there 
%is no doubt
can be no doubt that
% can be no question 
%%%%%%%% that 
the Second
Incompleteness Theorem is  fully
robust
% result 
from a purist 
%pristine 
mathematical perspective.
Yet,
it is still problematic to fully
% 
%  simultaneously
% % at the same time,
% it is 
% hard to 
% entirely
% 
dismiss
 Hilbert's 1926
suggestion that 
 some 
specialized forms of logics should
%declaration 
%% 
%% concerns
%% in $\,*\,$ 
%% that 
%% {\it ``the honor of human understanding''}
%% requires
%% examining
%% % explaining 
%% % considering
%% how  logic systems can
%% 
possess
a type of well-defined 
 knowledge about their
own 
internal
consistency.
(This is because it is
%  highly
very
 awkward to explain how
% and why
human beings 
are able to
%can 
%manage to 
motivate 
their 
own 
%cogitations,
cognitive process, otherwise.)

Thus the  distinguishing
between
the
 questions Q-1 and Q-2, combined with
 Example \ref{ex-2.5}'s ``SelfRef''formalism and
Proposition \ref{th-3.3}'s 
% employed 
``$\theta$ primitive'',
will usher in our new approach.
%and \ref{th-7.1} are 
%% 
%% bbbbb
%% 
\end{remm}

Before  discussing
our new ``IQFS'' methodology,
% Double-Formatted Logics, 
it should be mentioned
that other unusual interpretations of the Second Incompleteness
Theorem have followed
from Gentzen's perspectives about
transfinite induction 
under his $\epsilon_0$ ordinal
\cite{Ge36,Ta87}, the 
%% 
%% 
%% explore
%% how \cite{wwapal}'s results for a Single-Formatted logic
%% can be revised
%% % with our new $~\zzthe~$ function
%% under a
%% 
%% Before 
%% broaching
%%  this topic it should be mentioned that
%% %0fascinating
%% other approaches to
%% %efforts to partially 
%%  the Second Incompleteness Theorem 
%% % do
%% have centered around
%% 
 Kreisel-Takeuti's    ``CFA''
system \cite{KT74}
and also
the {\it interpretational frameworks} of
Friedman,
Nelson, Pudl\'{a}k and Visser
\cite{Fr79b,Ne86,Pu85,Vi5}.
These systems are unrelated to 
our 
%% main
%\cite{ww93}--\cite{ww14}'s 
methods.
%approach.
They
do not use
Kleene-like {\it ``I am consistent''} axiom-sentences.
Also,
they
%apply to 
employ
``cut-free'' logics
(rather 
than
a preferable Hilbert-style deductive apparatus).
%that 
%%%%%%%%%%% explored 
%%%%%%%%%%% in 
%%%%%%%%%%% \textsection \ref{ss32} ).
%%%we are considering).
%%
%%Instead, CFA uses the 
%%special
%%properties of ``second order'' generalizations of Gentzen's
%%{\it cut-free}
%%Sequent Calculus, 
%%and 
%%the
%%interpretational approach
%%formalizes how some systems 
%%recognize their
%% Herbrand consistency 
%%on localized sets of integers,
%%which 
%%unbeknownst to 
%%themselves,
%%includes all
%%integers.
%%
%%%These
%   alternate 
%%%approaches 
Their 
%alternate 
% very
 fascinating
perspective  
should
% certainly,
 be examined by researchers
interested in the
Second 
Incompleteness Theorem,
although 
%but
it is
%% 
% they are 
unrelated to 
our approach in the next section.

% \
% \particular
% \% the next section's 
% \%specific analysis of
% \type of
% \Hilbert-styled self-justifying effects,
% \studied in the current article.

%% systems 
%% formalizing
%% %verifying 
%% their
%% own consistency
%% %%%%%Definition \ref{def-2.2}'s
%% %%% approximate 
%% under
%% Hilbert-styled 
%% deduction.

%deduction.
%  Hilbert deduction.

%methods.
%formalism.

%% It is,
%% % They
%% %are, 
%% however,  not germane to the next section's
%% perspective.

%methodology.
%main formalisms.
%methods.
%results.

 % \baselineskip = 1.8 \normalbaselineskip 

%\section{
%\small 
%Improving \cite{wwapal}'s Results with a 
%``Double-Formatted'' Logic  }

% \section{Formalizing the First Half of our 2-Part Conjecture}
% the Statement     of the
\section{The ISCE and IQFS
Axiomatic Formalisms}

% \label{ss32}
\label{ss5}

The only aspect of our prior research that will be directly
related to our new IQFS axiomatic framework is
\cite{wwapal}'s ISCE formalism.  
The next several paragraphs will review 
\cite{wwapal}'s defintions of
ISCE, given in its Sections 3 \& 4.
% of \cite{wwapal}.
After reviewing ISCE's properties,
the remainder of this section will explain,
{\it intuitively,} how 
%the 
our new
IQFS framework
% ISCE 
can
be incrementally defined, by suitably
revising 
ISCE's 
prior
% formal 
%structure.
definition. (A subsequent {\it  formal proof,} 
affirming IQFS's 
% very
exact
self-justification properties, will
appear in 
 \textsection \ref{nn6}.)

% ISCE's formalism in a straightforward manner,

%% 
%% This section will
%% review \cite{wwapal}'s results in sufficient detail
%% so that a reader need not examine \cite{wwapal}'s formal 
%% text,
%% 
%% %%%%%definition of the ISCE axiom system. 
%% 
%% During our discussion, 
%% 

During our 
discussion,
%review of \cite{wwapal}'s results,
$~L^G~$ will
% once 
again denote
our Grounding-level language built out of our six basic
non-growth functions
% consisting  of 
(e.g. the
Subtraction, Division,
Maximum, Logarithm, Root and Count operations).
Also, $\,C_0\,$, $\,C_1\,$ and $\,C_2\,$
 will again denote
 %    again denotes
three constant symbols designating the
integers values of 
``0'', ``1'' and  ``2''.
In a context where Pred$(x)$ is an abbreviation for
``$\,x \,- \, 1\,$''
(or more precisely  
``$\,x \,- \, C_1\,$'' ), 
the ISCE axiom system from \cite{wwapal}
had
used
\eq{start}'s axiom
 statement 
to define
 $\,C_0\,$, $\,C_1\,$ and $\,C_2~$:
% these three constants:  
\begin{equation}
\label{start}
\mbox{Pred}(    C_0    )  =  C_0~ \, \wedge ~ \,
C_1 \neq C_0~ \, \wedge ~ \,
\mbox{Pred}(    C_1    )  =  C_0 ~ \, \wedge ~ \,
\mbox{Pred}(    C_2    )  =  C_1 ~ \, \wedge ~ \,
\forall \, x ~ \neg ~(~C_0 < x < C_1~) 
\end{equation}
%Also,
The  challenge 
\cite{wwapal} 
faced was its formalism could
not use any of the
%  conventional 
function-operations of
successor, addition or multiplication to infer the existence
of larger integers from the initial constants of 
$\,C_0\,$, $\,C_1\,$ and $\,C_2\,$. This was because
the Pudl\'{a}k-Solovay result $++$ 
indicated
  the 
presumption  successor is a total function 
precludes   
most
%axiom 
systems
from recognizing their
own Hilbert
consistency.

Our article
\cite{wwapal} 
considered two alternatives
%to a conventional Successor 
%function symbol 
for overcoming these difficulties,
 called
the  {\bf Additive} and  {\bf Multiplicative Naming}
conventions.
They defined 
some
further constant symbols $~C_3,~C_4,~ C_5,~ ...~$
and  $~C^*_3,~ C^*_4,~ C^*_5,~ ...~$
where 
%respectively
$~C_j~=~2^{j-1}~$ and $~C^*_j~=~2^{\,  2^{ \,j-2}}~$
under these two conventions.

The definition of these
% new 
constants
%  symbols
is 
easy
%straightforward
under $L^G\,$'s
% Grounding-level 
language.
% called $L^G~,~$
%all of whose function objects are non-growth primitives. 
This is because
Lines
\eq{newadd}  and \eq{newmult}
%had 
specify how
% that
two 3-way predicates, called
Add$(x,y,z)$  and  Mult$(x,y,z)\,$, 
%can 
%do
encode the identities of
% can be encoded to specify,  respectively,
$x=y+z$ and $x*y=z$.
Our additive and multiplicative 
% naming 
conventions
can,
%will,
 then, define 
 $~C_3,~C_4,~ C_5,~ ...~$
and  $~C^*_3,~ C^*_4,~ C^*_5,~ ...~$
%by  using
via
an infinite number of instances  of
%utilizing respectively
 \eq{addcov} and
\eq{multcov}'s
particular
axiom schemas, 
% defined below:
respectively:

{
%\small
\beq
\label{addcov}
\mbox{Add}(~C_{j-1}~,~C_{j-1}~,~C_{j}~)
\enq
\beq
\label{multcov}
\mbox{Mult}(~C^*_{j-1}~,~C^*_{j-1}~,~C^*_{j}~)
\enq}
The methodology in
 \cite{wwapal} 
%% employed \eq{addcov}  and  \eq{multcov}'s schema in a context where it 
presumed
% assured
%the Y of 
the ``names'' for its constants $ C_j $ 
and $ C^*_j $ 
had nice compact encodings using $O(~Log(j)~)$ bits.  
Its formalism calculated
%, thereby, 
the values of ``unnamed'' integers from 
named entities via the {\it non-growth} Subtraction and
Division primitives. For instance since $~20~=~32-8-4~,~$
the quantity 20 
can be encoded as $~C_6-C_4-C_3$.
%%%%%%%%%%%%%%%% under \eq{addcov}'s naming convention.

%% required
%% $O(~Log(j)~)$ bits.  
%% Thus, the length of these encodings was 
%% much 
%% smaller
%% than the respective
%% numbers
%% % magnitudes of 
%%  $2^{j-1}~$ and $2^{2^{j-2}}$ 
%% %that 
%% these constants represent.

\bigskip

The challenge \cite{wwapal} 
faced was to determine whether 
%it was possible to formulate 
self-justification
was possible
under 
either
\eq{addcov}'s
% ``Additive'' 
or  \eq{multcov}'s 
% ``Multiplicative'' 
%% naming
schema.
It found 
%that 
 \eq{multcov}'s 
multiplicative
% naming 
convention was  incompatible 
with self-justification (due to its 
%%very 
speedy growth rate),
but
%In contrast,
\eq{addcov}'s additive 
% naming 
schema did,
conveniently,
 permit self-justification. 

\bigskip

Our new proposed Double-Formatted form of a self-justifying
axiom system is easiest to describe, if we first
review \cite{wwapal}'s Single-Formatted formalism
and then 
%incrementally 
refine it.

\bigskip

The extension of our base-language $~L^G~$
that includes the Additive Naming Convention (ANC)'s
additional constants 
 $~C_3,~C_4,~ C_5,~ ...~$ will be called
an {\bf ANC-Based Language} and be denoted
as  $~L^{ANC}~$. 
Also if
 $\, t \,$ denotes   any term in $\, L^{ANC} \,$'s   
language, then 
the quantifiers in 
the two wffs of
$~ \forall ~ v \leq t~~ \Psi (v)~$ and
$\exists ~ v \leq t~~ \Psi (v)$
will be  called $\, L^{ANC} \,$'s   
{\bf ``Bounded 
Quantifiers''}.

\medskip

\begin{deff}
\label{def-3.8}
\rm
The analogs of a conventional arithmetic's
$\Delta_0$, $\Pi_n$ and $\Sigma_n$ 
formulae
in the
language $L^{ANC}$ will be denoted as
$\Delta^{ANC}_0$,
 $\Pi^{ANC}_n$
 and $\Sigma^{ANC}_n$.
Thus,
a formula is defined to be
$\Delta^{ANC}_0$  iff all its quantifiers 
%are so  bounded.
satisfy the prior paragraph's bounding condition.  
The
%%%%%%%%%  below 
definitions 
of  $\Pi^{ANC}_n$ and $\Sigma^{ANC}_n$
formulae, 
%% specified below 
in Items 1-3, are also 
 quite 
conventional:
\bee
%\small
%\parskip -2 pt
%\baselineskip = 0.8 \normalbaselineskip 
\item
Every  
$\Delta_0^{ANC}$ formula is considered to
be 
also 
a
$\Pi_0^{ANC}$  and 
an
$\Sigma_0^{ANC}  $ expression.
%% 
%% ``$~\Pi_0^{ANC}~ \,$''  and 
%% %  also 
%%  ``$~\Sigma_0^{ANC}~ \, $''. 
%% 
\item
A
formula
is  called
 $ \,\Pi_n^{ANC} \,$
when it
% is 
can be
encoded as 
$\forall v_1 ~ ...~ \forall v_k ~ \Phi$  
where
%with
$\Phi$ is  $\Sigma_{n-1}^{ANC}$
\item
A formula
is  called
 $\Sigma_n^{ANC}$
when it can be encoded as 
$\exists v_1~ ...~ \exists v_k ~ \Phi,$  where
$\Phi$ is  $\Pi_{n-1}^{ANC}$.
\ene
\end{deff}

%%\begin{deff}
%%\label{def3.9}
%%\rm

\bigskip

Given an initial axiom system $\beta,$
the Theorem 3 of \cite{wwapal}
% defined 
%formalized
did formalize
a 
self-justifying logic, called
{\bf ISCE$(\beta)$}, 
that could prove all 
$~\beta\,$'s $\Pi_1^{ANC}$ theorems and 
which could also
verify its own consistency under any Hilbert-style deductive
apparatus
(including the $d_E$, $d_M$ and $d_H$ deductive methodologies) 
The axiom basis for ISCE$(\beta)$ 
was comprised, formally,  of the following four
distinct
groups of axioms:
%
%  \newpage
\begin{description}
 \parskip 0pt
\item
{\bf GROUP-ZERO:} 
This 
schema
%  axiom group 
will
use \el{start}'s axiom to define the constants of
$\,C_0\,$, $\,C_1\,$ and $\,C_2\,$ 
and 
%employed
%an infinite number of instances of
\el{addcov}'s Additive Naming 
schema
%convention
to define 
 the further constants
 $    C_3,  C_4,  C_5,  ...   $  
\item
{\bf GROUP-1:}
It is convenient  to
define
 ISCE's Group-1
% and Group-2
% axiom
schema
 using a notation that 
will support \cite{wwapal}'s Theorem 3
in a 
% slightly
 more 
general sense than
had
appeared in \cite{wwapal},
%%% under a slightly different notation convention,
% is transparently equivalent
% (but slightly different) from \cite{wwapal}'s counterpart,
so that
% a 
our
% the
new
% proposed
%%%% second
 ``IQFS''
formalism
(appearing later in this section)
% shall
will
%proposal shall
% framework will 
be easier to describe.
Let us
therefore
 say  a $\Pi_1^{ANC} $  sentence is {\bf Simple}
iff the only built-in constants it employs are
$\,C_0\,$, $\,C_1\,$ and $\,C_2$.
Then ISCE's Group-1 scheme will be
allowed to
 be any finite set of 
simple $\Pi_1^{ANC} $  axioms, called $~S~,~$
that is consistent with 
the
Group-zero schema and
which
 has
the following properties:
\bee
\item
  The union of $~S~$ with ISCE's Group-Zero
axioms
will be capable of proving all $\Delta^{ANC}_0 $  
statements which
are true.
\item
  The union of $~S~$ with ISCE's Group-Zero
scheme
will also be capable of proving 
that
the ``=" and ``$\leq$" predicates
% own 
support
their conventional
transitivity, reflexivity, symmetry and total ordering
properties.
\ene
Any finite set of 
$\Pi_1^{ANC} $  axioms
with the above properties can be used to define $~S~$
and 
support
%prove 
an analog of
\cite{wwapal}'s Theorem 3,
by 
%a trivial 
an easy
generalization
\footnote{A 
%formal 
proof of this 
specific
generalization of
\cite{wwapal}'s results is 
%absolutely 
entirely
routine.
%  and omitted here for the sake of brevity.}
For the sake of brevity, it  is
omitted.} 
%
%here,}
 of
the methodologies from Sections 3 and 4 of
\cite{wwapal}. (Thus,
any such finite set $~S~$ supporting Conditions (1) and (2)
can  be employed
%by
to define
ISCE's
Group-1 axiom schema.)
%%% 
%%% and it is unimportant which
%%% particular defining
%%%  set is used.
% 
% BBB111
%  
% This
% % schema
% axiom group 
% consisted of a finite
% set of 
% $\Pi_1^{ANC} $ axioms
% %, CALLED $F$,
%  defining ISCE's
% Grounding  function primitives. 
% %This means that
% For each  such function $G$ and set of numbers 
% $    {k},   {k_1},    {k_2}, ...    {k_m}$,
% %the combination of 
% the Group-Zero and Group-1 axioms 
% %must
% will
%  imply 
% $ G(    {k_1},    {k_2}, ...    {k_m}) \,=\,    {k}   $ when
% this sentence is true
% \footnote{ \f55 
% Our
%  $\Pi^{ANC}_1$
% encoding for the
% Group-1 scheme  needs,
% technically, 
% % employ
%  only
% employ
%  the three constant symbol $C_0$, $C_1$ and $C_2$ for the
% union of all 
% the
%  Group-Zero and Group-1 axioms
% to satisfy 
% their
% %its
% %the 
% above requirements.} .
% The Group-1 schema
% of \cite{wwapal} 
% will also
% assign the ``=" and ``$<$" predicates
% their conventional
% %  logical
% properties.
% %footnoted property.}
% %%
% %%(Any finite 
% %%set of  $\Pi_1^{ANC} $
% %%sentences  meeting these conditions is
% %%suitable.)
% %%
\item
{\bf GROUP-2:}
Let
$\ulcorner \, \Phi \, \urcorner$ denote $\Phi$'s G\"{o}del number, and
$\mbox{HilbPrf}_{ \beta  }(x,y)$
denote a 
%%%%%%%%%%%  $\Delta _0^{ANC+}$ 
$\Delta _0^{ANC}$ 
formula indicating  $y$ is a
Hilbert-styled
proof
from axiom system $\beta  $ of the theorem 
$x$.
%
% Suppose that
%$~\beta~$ uses the same Grounding function symbols as
%ISCE$^{ANC}(\beta)$,
%and it therefore generates
%a set of 
%% $\Pi_1^{ANC+} $ theorems.
% $\Pi_1^{ANC} $ 
%theorems.
%
For each 
%$\Pi_1^{ANC+} $
$\Pi_1^{ANC} $
 sentence  $\Phi$,
the Group-2 schema
for ISCE$(\beta)$ 
%
%was defined in  \cite{wwapal} 
%did 
will 
contain
% an 
%one
a $\Pi_1^{ANC} $
axiom of the form:
%% 
%% \begin{equation}
%% \small
%% \label{group2nold}
%% \forall ~x~\forall ~y~
%% ~~\{~~[~~ \sigma_{~ \ulcorner \, \Phi \, \urcorner
%%  ~}(x)~\wedge~
%% \{~ \mbox{HilbPrf}~_\beta
%% ~(~ x ~,~y~)~~]~~
%% \Rightarrow ~~ \Phi~~ \}
%% \end{equation}
%% % {\bf IMPORTANT CLARIFICATION:} 
%% %{\small 
%% %%{{\bf DECIPHERING LINE \eq{group2nold}:$~$}
%% {{\bf Clarification:$~$ }
%%  \el{group2nold} is {\it helpful}
%% because   ISCE(\beta)$  can             infer
%% \eq{group2old}'s {\it simpler statement}
%% directly
%% from the combination of
%% \eq{group2nold},
%% % it,
%% the Group-1 schema and \el{deltf}'s definition of
%% ``$~\sigma~$''.}
%% 
\begin{equation}
% \small
\label{group2old}
\forall ~y~~~\{~ \mbox{HilbPrf}~_\beta
~(~ \ulcorner \Phi \urcorner ~,~y~)~~
\Rightarrow ~~ \Phi~~\}
\end{equation}
\item
{\bf GROUP-3:}
This last part of
%%%%%%%%%%%%%%%  
\cite{wwapal}'s
ISCE$(\beta)$
formalism
is
% was
 a single 
self-referencing
$\Pi_1^{ANC}$
sentence  
stating:
 %% essentially declaring:
\begin{quote}
% \small
%%%%%%%%%%%%% $ \oplus ~ \oplus ~~~$
$ \oplus  \oplus ~~~$
 ``There 
%is 
exists
no
Hilbert-style proof of 0=1 from the union of the Group-0, 1 and 2
axioms  with {\it THIS  SENTENCE} (referring to itself)''.
\end{quote}
\end{description}
%{\bf CLARIFICATION:}
{\bf Clarifying $ \oplus  \oplus$'s Meaning:}
 $~$Several of our articles
\cite{ww1,ww5,wwapal,ww9}
employed
self-referential
  $\Pi_1^{ANC}$ constructions,
similar 
to 
 $ \oplus  \oplus \,~$,
as the Example \ref{ex-2.5} 
%had mentioned.
did mention.
A reader can, thus,  find
several
%detailed
% slightly different
 illustrations about how 
$~ \oplus  \oplus ~ $
%  $\, \oplus  \oplus $'s 
%   self-referential statement 
is encoded in  these articles.

% 
% Each of these articles provide examples of
% how analogs for
% $\, \oplus  \oplus $'s 
% self-referential
% statement
% are encoded.

% If the reader wishes to see
% a formal encoding  for
% $\, \oplus  \oplus $'s 
% %self-referential
% % Fixed-Point 
% statement,
% %it 
% one such example
% is provided by  
% \cite{wwapal}'s 
%  Lemma 1.
% 

\begin{deff}
\label{def-3.9x10}
\rm
Let $~I(~\bullet~)~$ denote
an operation that maps
an initial axiom basis $\, \beta \,$ onto an alternate
system  $\,I(\beta)\, $.
(One example of
such an operation is the
  ISCE$( \, \bullet \, )$ 
framework,
that maps 
an initial axiom basis of 
  $~\beta~$ onto 
the alternate formalism of
 ISCE$(\beta).~)~$ 
Such an operation  $~I(~\bullet~)~$
is  called {\bf Consistency Preserving}
iff  $\,I(\beta)\, $ is consistent whenever 
the union of
 $\beta$ with the Groups 0 and 1 axiom schemas is
consistent.
\end{deff}

%Most of our research in 
% \cite{ww93}-\cite{ww14}
% has 

\begin{exx}
\label{ex-5.3}
%\label{ex-basis}
\rm
Several of our research projects 
%centered around
had employed
 \dfx{def-3.9x10}'s
framework.
For instance, 
%% 
%% the 
%% 
%% 
%% Its
%% %%%  main
%% % central 
%% focus in
\cite{wwapal} 
demonstrated
%consisted of showing
 the   ISCE$( \, \bullet \, )$ 
mapping was consistency preserving.
Thus if PA+ denotes the extension of
Peano Arithmetic that 
includes
PA's traditional  Addition and Multiplication
functions
%% 
%% 1n addition to the conventional
%% functions of addition and multiplication
%% contains 
%% 
%% 
plus $L^G\,$'s six
added
%previously mentions
 Grounding-level function 
primitives,
%functions,
then
 ISCE$( \, $PA+$ \, )$
must
%will 
% be 
automatically
be
 consistent
(because  PA+ was consistent).
% consistent whenever PA+ is consistent.
Hence while Peano Arithmetic is unable to
verify its own consistency (on account of G\"{o}del's
seminal 1931 discovery), it is 
known to be
sufficiently agile to
prove the following relative-consistency statement:
\begin{center}
%% \small
$\#~~~$ If PA is consistent then 
 ISCE$( \, $PA+$ \, )$ is  
 self-justifying.
 \end{center}

%This
The
preceding
% above
% statement
 relative-consistency statement
%does offer 
provides
a partial 
positive
answer to
the
Q-2 version of Hilbert's Second  Question.
%This is because it formalizes one respect 
It 
captures
% Brad change encapsulizes 
one
% positive 
respect
in which
%such as  
ISCE$( \, $PA+$ \, )$
can {\it appreciate} its own consistency.
This respect is, obviously,
% only
of a 
sharply
limited nature
because $++$'s generalization of the
% Second
Incompleteness Theorem indicates 
%that
no Type-S arithmetic
can
% simultaneously 
recognize 
% {\it both} 
its Hilbert consistency and
view
%take
successor 
as
%to be
 a total function.
%The consistency-preservation property of 
% ISCE$( \, \bullet \, )$
%dies, however,
It does,  however,
 raise 
an
%the following
enticing
 question:
\begin{quote}
$\# \, \#~ $
Can the infinite number of
distinct
 constant symbols, employed by
ISCE's Group-Zero schema, be reduced to a finite size
by a Type-NS Self-Justifying Logic
in some type of formal manner?
\end{quote}
The
next two sections of this paper
 will outline how 
%an 
a quite
encouraging
answer to 
$\, \# \, \# \, $'s query
shall
%will
%is likely to
%%%should,
% conveniently
arrive,
%be plausible
when one 
% carefully
%delicately 
modifies ISCE's formalism
with  the Q-function operative of $~\zzthe~$.
\end{exx}

%\smallskip

\begin{deff}
\label{def-3.10}
\rm
Let $L^Q$ 
 once 
again denote the 
special
extension of
$~L^G\,$'s Grounding language that includes
the
added new
% further
 Q-function symbol of $\, \theta $.
Then
$\Delta^Q_0$,
 $\Pi^Q_n$ and $\Sigma^Q_n$
will
%  , intuitively,
%similarly
 denote the
%  1-to-1 
analogs of
\dfx{def-3.8}'s
$\Delta^{ANC}_0$,
 $\Pi^{ANC}_n$ and $\Sigma^{ANC}_n$'s
formulae
%in
under 
$~L^Q\,$'s 
modified
 language.
In particular, if $~\Phi~$ 
is one of an
$\Delta^{ANC}_0$,
 $\Pi^{ANC}_n$ or $\Sigma^{ANC}_n$
formula,
then
% the formula 
$~\Phi^Q~$
will be called
% respectively
$\Delta^Q_0$,
 $\Pi^Q_n$ or $\Sigma^Q_n$
when
%if 
it
differs from $~\Phi~$ 
%  only
by
exactly
 replacing each constant $~C_J~$
from the set $~C_3,C_4,C_5...~$ 
with Line \eq{ej-def}'s 
% mathematically equivalent term of
term  $~E_{J-1}~$. 
\end{deff}

\parskip 2pt

%% 444444444444444

\begin{example}
\label{ex-3.11}
\rm
Suppose $~\Phi$ 
is one of a 
$\Delta^{ANC}_0$,
 $\Pi^{ANC}_n$ or $\Sigma^{ANC}_n$
sentence that employs the three constant symbols
of $C_4$,   $C_6$ and  $C_{10}\,$
for
 representing the
three numbers
of 8, 32 and 512. 
Let us recall 
that
 $E_3$,   $E_5$ and  $E_9\,$
%  do
formulate these three quantities
under  Line \eq{ej-def}'s notation.
Then $~\Phi^Q$  will have an 
identical definition as
 $~\Phi$ 
except each $C_j$ is replaced by
$E_{j-1}$.

\smallskip

A formula is,
moreover,
 defined to lie in one
of the 
$\Delta^Q_0$,
 $\Pi^Q_n$ or $\Sigma^Q_n$
classes 
{\it if and only if} it is constructed in
exactly
 such a manner.
This fact
% brad assures 
ensures
that all the
main 
primary
terms employed in our
three 
major
classes of sentences are 
{\it ``Observable''} terms.
Hence ``Unobservable'' ground terms are allowed in
$~L^Q\,$'s language, but {\it they are excluded}
from occurring as the 
{\it ``end-product 
primary''}
terms in the 
$\Delta^Q_0$,
 $\Pi^Q_n$ or $\Sigma^Q_n$
theorems 
that 
%%will now be discussed.
%it proves !
do
encapsulate
%formalize
 the {\it intended use of 
%its
our
formalism.}
We will need one more preliminary definition before 
%%our new 
``IQFS''
% axiom system
% and its language $L^Q$ 
can be
% formally 
defined.
\end{example}

\smallskip

\begin{deff}
\label{def-pain}
\rm
Let us recall that Example \ref{ex-2.1} defined three
examples of Hilbert-style deductive methodologies.. 
The symbol $d_{ER}~$ 
will denote the particular version of
Hilbert-style deduction that we will use in  this article.
It will be a minor
revision
of Enderton's $d_E$ apparatus
(that was 
mentioned
% formalized 
in Example \ref{ex-2.1}).
 Unlike Enderton's schema,
it will not take all tautologies as axioms.
Instead, it will use the following
three
 axiom schemas, from
Mendelson's textbook \cite{Mend},
to prove all tautologies as
theorems:   
\bed
\item[   I.] $~B~ \Rightarrow ~ ( ~C ~ \Rightarrow ~ B~)~$
\item[  II.] $[~B~ \Rightarrow ~ ( ~C ~ \Rightarrow ~ D~] ~ \Rightarrow ~
[~ ( ~B ~ \Rightarrow ~ C~)  ~ \Rightarrow ~
( ~B ~ \Rightarrow ~ D~) ~]$
 \item[ III.]  $[~(~ \neg C~) \Rightarrow ~(~ \neg~ B~) ~]  ~ \Rightarrow \,$
$ \{ ~ [~(~  ~\neg~ C ) ~  ~ \Rightarrow ~ B ]  ~ \Rightarrow ~ C~~\}$

\ennd
This approach is preferable because the problem of identifying
tautologies is NP-hard (e.g. 
it is preferable to 
avoid viewing tautologies as axioms and to instead employ
the axiom schemas I-III to 
prove tautologies as theorems).
Thus, the logical axioms of 
$d_{ER}~$  shall consist of the axiom schemas
I-III plus the below schemas  2-4 
from \textsection
 2.4 of  Enderton's textbook.
The deductive apparatus $d_{ER}~$  will also
follow \cite{End}'s example
by treating $\forall~ x ~\Psi~$ automatically as a logical axiom
whenever $~\Psi~$ is a logical axiom.
\bed
\item[ 2. ] $~\forall ~x~ \phi(x)~  \Rightarrow ~  \phi(t)~ $
when $~t~$ is substitutable for $~x~$ in $\phi$.
\item[ 3. ] $~\forall ~ x~ (~\phi ~ \Rightarrow ~ \psi ~) ~~ 
\Rightarrow ~~~\{~ (~ \forall ~x~ \phi ~)~ 
 \Rightarrow ~(~ \forall ~x~ \psi ~)~ \}$ 
\item[ 4. ] $~\psi~  \Rightarrow ~ \forall x ~ \psi $ when $x$
does not occur free in $\psi$. 
\ennd
The symbol  $~L_{d_{ER}}~$ 
will henceforth denote this  set  of six 
basic
logical axiom schemas
(e.g. items I-III combined with 2-4).
Our
${d_{ER}}$ 
deduction methodology
 will follow Enderton's example
by using  modus ponens as
% the 
its
only 
employed
 rule of inference
\footnote{Unlike several other 
textbooks, Enderton \cite{End} does not view
generalization as a
formalized rule of inference
 because he takes Items (3) and (4)
as logical axioms.}.
\end{deff}

\smallskip 

The remainder of this section will 
combine the formalisms of Defintions \ref{def-3.10} and 
\ref{def-pain} to define our new axiom system, 
 IQFS($ \beta$),  
that can prove all $\, \beta \,$'s
$\Pi_1^Q$ theorems as well as  recognize its own consistency
under ${d_{ER}}$'s form of Hilbert-style deduction.
In essence,  IQFS($ \beta$)$~$  will be the direct
analog of \cite{wwapal}'s  ISCE($ \beta$) axiom system
(when the logic's base language is changed from 
   $\,L^{ANC}\,$ to $L^Q$ and our focus changes from 
$\Pi_1^{ANC}$ to $\Pi_1^Q$ theorems). 

\smallskip 

\begin{deff}
\label{def-3.12}
\rm
The acronym ``IQFS'' stands for 
``Introspective Q-Function Semantics''.
In a context where $~\beta~$ is 
%some i
an initial axiom
system that proves theorems 
%under
in 
the 
language  $L^Q$, the 
system,
%formalism 
 IQFS($\, \beta \, )$
%$~$
will 
be designed to be a
 4-part 
formalism,
analogous to  ISCE($\beta$),
that can prove all $\beta$'s
$\Pi_1^Q$ theorems and recognize its own consistency.
In essence,  IQFS($ \beta  )$'s definition will be identical
to
% that of  
ISCE($\beta$), except for the following
four changes (three of which are 
%quite 
trivial):
\bed
\item[  a.  ]
The 
{\bf 
GROUP-ZERO} schema of
 IQFS will
differ from ISCE's analog
by replacing 
\el{addcov}'s ``Additive Naming'' schema with  
the
Up-Walking axioms,
%given in
from
 Lines  \eq{walk1}--\eq{walk4}.
(This is because
the language  $L^Q$ differs from
 $L^{ANC}$ by 
having the 
 Q-function operator of $~\zzthe~$
define the formal quantities that are represented by
the constant symbols
of $~C_3,C_4,C_5~~....~$ 
under  $L^{ANC}.~~)$ 
 In order to assure that the 
four
Up-Walking axioms
are adequate to determine $E_j$ 
is numerically equivalent to  $C_{j+1}$, our new Group-Zero
schema will include
% once again 
\el{start}'s 
axiom (for defining the starting
 constants of
$\,C_0\,$, $\,C_1\,$ and $\,C_2\,$) 
as well as \eq{halfax}'s added axiom:
% sentence. 
\beq
\label{halfax} 
%\forall ~x\, \geq \,2 ~~ \frac{x}{ \mbox{Half}(x)}  ~=~2
\forall ~x  ~~ \{ ~~[~\mbox{Power}(x) ~\wedge ~ x \geq 2~]~~\Rightarrow~~\frac{x}{ \mbox{Half}(x)}  ~=~2~\}
\enq  
%% 
%% Otherwise both
%% these
%% Group-Zero 
%% schemes will be 
%% identical.
%% Thus,
%%  they 
%% will 
%% both 
%% use \el{start}'s axiom to define the 
%% three initial constants of
%% $\,C_0\,$, $\,C_1\,$ and $\,C_2\,~$.  
%% 
\item[  b.  ]
The {\bf GROUP-1} scheme of IQFS will be {\it identical} to ISCE's counterpart,
except it will reflect Item a's modification of the Group-Zero
scheme for $\, L^Q\,$'s 
particular
revised language
(e.g. 
the footnote \footnote{\f55
In a context where a $\Pi_1^Q $  sentence is
called  ``Simple''
when it contains no $~E_j~$ term with $~j \geq 2~,~$ 
the Group-1 scheme of IQFS will be analogous to
ISCE's counterpart by consisting of
 any finite set of 
simple $\Pi_1^Q $  axioms, called $~S^*~,~$
that is consistent with Group-zero schema and
which
 has
the following properties:
\bee
\item
  The union of $~S^*~$ with IQFS's Group-Zero
axioms
will be capable of proving all $\Delta^Q_0 $ 
statements which  are true.
\item
  The union of $~S^*~$ with IQFS's Group-Zero
scheme
will also be capable of proving 
that
the ``='' and
 ``$\leq$'' predicates
support their conventional
transitivity, reflexivity, symmetry and total ordering
properties.
\ene
The above two properties are the 1-to-1 analogs of
their counterparts used by ISCE's Group-1 scheme.
As was the case with ISCE's formalism,
any finite set
of simple 
$\Pi_1^Q $  axioms
with the above properties can be used to define $~S^*~.~$
Once again,
it is unimportant which 
particular
finite-sized
definition for $S^*$
% such set 
is used.} describes
this  
easily
%quite
 straightforward 
% resulting
%% quite
revision of
ISCE's Group-1 scheme.)
\item[  c.  ]
All the $\Pi_1^Q$ axioms lying in IQFS's
%Group-1 and 
{\bf GROUP-2} scheme will be identical to their counterparts
under ISCE, except  they
will
 employ 
\dfx{def-3.10}'s machinery for translating 
 $ \,\Pi_1^{ANC} \,$ 
sentences into their equivalent
 $ \,\Pi_1^Q \,$ counterparts.
In particular, let us assume  the
G\"{o}del number ``$~ \ulcorner \Phi \urcorner ~$'' is defined via
the ``Byte-Style'' encoding method (defined in the next
section)
and 
$\mbox{HilbPrf}_{ \beta  }(x,y)$
denotes a 
%%%%%%%%%%%  $\Delta _0^{ANC+}$ 
$\Delta _0^Q$ 
formula indicating  $y$ is a
Hilbert-styled
proof of the theorem 
$x$
from axiom system $\beta  $. 
Then for each 
$\Pi_1^{Q} $
 sentence  $\Phi$, 
$~$IQFS$(\beta)$'s  Group-2 schema
will 
contain
 a $\Pi_1^{Q} $
axiom of the form:
\begin{equation}
\label{group2iq}
\forall ~y~~~\{~ \mbox{HilbPrf}~_\beta
~(~ \ulcorner \Phi \urcorner ~,~y~)~~
\Rightarrow ~~ \Phi~~\}
\end{equation}
\item[  d.  ]
The {\bf GROUP-3} axiom of  IQFS
will be 
essentially
similar to ISCE's Group-3 
{\it ``I am consistent''} 
% axiom-statement, 
axiom,
except 
the latter's notion of ``I'' will 
trivially
reflect the above
changes in the Groups 0, 1 and 2 schemes. 
It
%Thus, the new
%Group-3 axiom 
will,
thus,
be a $\Pi_1^Q$ sentence declaring that:
\begin{quote}
%ggg333
$\oplus~\oplus~\oplus~$
{\it ``Under Definition \ref{def-pain}'s Hilbert-style deductive methodology,
there will exist no
proof of 0=1 from the union of the preceding axioms
with THIS SENTENCE (looking at itself)''.}
\end{quote}
Readers  familiar with our previous
papers \cite{ww93,ww1,ww5,wwapal}
 about 
{\it ``I am consistent''} 
% axiom-statement, 
axioms  
should be able to 
 appreciate,  intuitively, 
it
 is feasible
to formulate  a $\Pi_1^Q$ 
encoding for
  $\oplus\, \oplus\, \oplus $.
% 's self-affirming declaration.
$~$This fourth part of our definition of the IQFS's axiom system
should, however, be articulated in 
%thorough 
greater detail because it is more complicated
than
 IQFS$(\beta)$'s
other three parts. 
Therefore, 
the next chapter's \lem{lem6.3}
will formalize
% , in detail,
%exactly 
precisely
 how
  $\, \oplus\, \oplus\, \oplus\, $
%can receive
is endowed with
 a $\Pi_1^Q$ 
encoding.
\ennd
\end{deff}

% DELETE NEXT SENTENCE ????
% For the benefit of
% those readers, less familiar with this topic,
% % subject,
% the next chapter's \lem{lem6.2} will explain how
%   $~\oplus~\oplus~\oplus~$
% can receive a formal  $\Pi_1^Q$ 
% encoding.

%This  article's 
Our
% underlying  
main
goal
%gist 
will be to show that
IQFS
satisfies
%does satisfy 
Definition \ref{def-3.9x10}'s 
Consistency-Preserving property, 
analogous to
\cite{wwapal}'s
%  prior  
ISCE methodology.
% from our earlier paper \cite{wwapal}.  
The 
particular
virtue of IQFS is that
it needs only three starting constants, $C_0$, $C_1$,  and $C_2$, to define
the full infinity of natural  numbers, 
unlike ISCE. 
(Moreover, our 
\phx{th-3.3} will  imply
IQFS can encode
every integer $\,n\,$ 
efficiently
with a   term $T_n$ 
%that has
with
  an
%% that has no more than an 
$O\{~[~$Log$(n)~]^3~\}$ length,
and \phx{th-7.1}
will show this quantity 
% $O\{~[~$Log$(n)~]^3~\}$ length 
%can be reduced 
can be reduced
to a {\it yet-more-efficient}
% essentially 
logarithmic size,
when ground terms are encoded as directed acyclic graphs.)

%% {\bf MAYBE REMOVE ABOVE parenthesis SENTENCE}

\section{Main Results about IQFS}
\label{nn6}

Our main results about
the IQFS
formalism will be proven in this chapter.
Its discussion will be divided into two
parts where:
\bee
\item Section \ref{6x.1}
will formalize the encoding of 
  $\, \oplus\, \oplus\, \oplus\, $'s Group-3 axiom
(e.g. 
the Item (d) from Definition \ref{def-3.12}).
\item Section \ref{6x.2} will display our
main  theorems about 
IQFS's consistency-preservation property.
\ene

\subsection{The  $\Pi_1^Q$  Encoding for IQFS's Group- 3 Axiom} 

\label{6x.1}

Before discussing how to endow  $ \, \oplus \, \oplus \, \oplus \, $ 
with a  $\Pi_1^Q$ encoding, it is necessary
to first outline
our methodology for generating a  proof's  
 G\"{o}del number. It
will be
% roughly
analogous to the natural B-adic encoding
methods
% methodology
used by  Buss \cite{Bu86},
H\'{a}jek-Pudl\'{a}k \cite{HP91} and
Wilkie-Paris \cite{WP87} --- insofar as
the
number of utilized bits to encode a
semantic object will be  approximately
proportional to the length of such an expression written
by hand. 
We will say a byte is a string of 6 bits, and encode each proof
as an integer, written in base 64, comprising a sequence of bytes.
Our encoding
 will use less than 32 logical symbols; thus, each symbol
will be encoded by an unique integer between 32 and 63.
These symbols shall
% will
 include:
\begin{enumerate}
\small \baselineskip = 1.1 \normalbaselineskip  \parskip    0pt
\item The standard connective symbols of
$\wedge ,~ \vee ,~ \neg ,~ \Rightarrow ,~ \forall$
and $~ \exists$.
\item The 
left and  right parenthesis symbols, and also a comma symbol (to separate
terms) and a period symbol (to specify the end of a 
sentence)
\item
Seven function symbols for representing
our six grounding functions 
and the new $\theta$ primitive.
\item 
The relation symbols of
``$~=~$'' and ``$~ \leq ~$''.
\item The symbol $~ \hat{V} ~$  for designating
the presence of a basic  variable symbol.
\item Three constant symbols to represent the built in constants
of 0, 1 and 2.
\end{enumerate}
The  $~ \hat{V} ~$  symbol 
is somewhat different from our other
logical symbols
because there are an infinite number of different variable  names
that may occur in a proof. The $j-$th variable will
be represented by a sequence   of $~1\, + \,\lceil \,$Log$_{32}\, (j+1) \,
\rceil\,$
bytes, where the first byte encodes the  $~ \hat{V} ~$  symbol and
the remaining bytes encode $~j~$ as a Base-32 number.

\begin{deff}
\label{def-6.1}
\rm
Let $\alpha$ denote a set of proper axioms.
Then the symbol   AxiomCheck$_\alpha(x)$ will 
denote
 a  predicate formula
that yields a value of True
when  $~x~$ is the
G\"{o}del number of one of $\alpha$'s axioms.
Also, ProofCheck$_\alpha(s)$  will denote a formula
that is True when $~s~$ represents a
byte string for a proof,
whose proper axioms come from $~\alpha~$ and whose
logical axioms and deductive apparatus were formalized
by  \dfx{def-pain}'s ``$d_{ER}$'' formalism.
The relationship between these two concepts is formalized
by Lemma  \ref{lem6.2}.
\end{deff}

% \bigskip

\bel
\label{lem6.2}
% \label{lem6.2}
The predicate  ProofCheck$_\alpha(s)$  is 
  $\Delta_0^Q$ encodable whenever
 AxiomCheck$_\alpha(x)$ is 
  $\Delta_0^Q$ encodable. {\rm  (The latter condition
does
 hold
for the majority of r.e. axiom systems $\alpha$.)}
\enl

{\bf Proof Summary:} It is well known that Lemma \ref{lem6.2}
holds for classic arithmetics that use 
  $\Delta_0$ encodings,
instead of   $\Delta_0^Q$ 
encodabilities. (There are many proofs of this fact.
One can  use the fact that Wrathall's notion of LinH functions
are known \cite{HP91,Kr95,Wr78} to have
the property that a formula has a  $\Delta_0$ encoding
if and only if 
it can be corroborated by
 a LinH decision procedure.)
A similar paradigm trivially generalizes for
 $\Delta_0^Q$ encodings. Other routine details
about \lem{lem6.2}'s proof
 are
omitted for the sake of brevity.
 $~~~\Box$  

\lem{lem6.2} 
is obviously trivial, but it was worth mentioning
because we will use it during one of the footnotes in
Lemma \ref{lem6.3}'s proof.

\bel
\label{lem6.3}
The Group-3 axiom
of $~$ IQFS$(\beta)$,
  introduced informally
by \dfx{def-3.12}'s statement
$\oplus \, \oplus \, \oplus  $,
can be formally encoded as a $\Pi_1^Q$ axiom  sentence.
\enl

{\bf Proof.}
Let
 UNION($\beta$)  denote the union of 
IQFS$(\beta)$'s Group-Zero,
Group-1 and Group-2 axioms
(given in Items a-c of \dfx{def-3.12}) ).
We will also use the  notation conventions below:
\begin{description}
\item{i. }   $\mbox{Prf}_{~UNION(\beta)}( \, t \, , \, p \,)$ 
will denote a formula  designating
that $~p~$ is  a proof of the theorem $~t~$ from the axiom
system  UNION($\beta)$, using 
\dfx{def-pain}'s deduction method $d_{ER}$.
\item{ii.}  
$\mbox{ExPrf}_{~UNION(\beta)}~( \, h \, , \, t \, , \, p \,)$  
will be a formula stating that
$~p~$ is a proof
(using deduction method $d_{ER}$)
 of
a theorem $~t~$ 
from the union 
of the axiom system
UNION$(\beta)$ with  an added axiom
%sentence 
whose G\"{o}del number equals
$\, h \,$. 
\item{iii.} 
 $\mbox{Subst} \, ( \, g \, , \, h \, )$ will denote
G\"{o}del's
classic substitution formula --- which yields TRUE when $\, g \,$
is an encoding of a formula
and $\, h \,$ is an encoding of a sentence
that replaces all occurrence of free variables in $\, g \,$ with
the G\"{o}del-encoded term $~\underx{g}~$.
\item{iv.} 
$\mbox{SubstPrf}_{~UNION(\beta)}~( \, g \, , \, t \, , \, p \,)$  
will denote the natural 
hybridizations of the constructs from Items (ii) and (iii).
It will yield a Boolean value of TRUE exactly when there
exists some integer $~h~$ 
simultaneously 
satisfying
{\it both}
the conditions
  $\mbox{Subst} \, ( \, g \, , \, h \, )$ and
$\mbox{ExPrf}_{~UNION(\beta)}~( \, h \, , \, t \, , \, p \,)$.

\end{description}
Each of (i)--(iii)
can be
easily
\footnote{ In particular, the $\Delta_0^Q$ 
encodings for Items (i) and (ii)
are immediate consequences of
\lem{lem6.2}.}
 encoded as  $\Delta_0^Q$ 
 formulae .
In such a context, \el{encode} illustrates
one possible $\Delta_0^Q$ 
encoding for 
$\mbox{SubstPrf}_{~UNION(\beta)} \,(   g   ,   t   ,   p  )$'s 
graph. (It is
equivalent to 
$~$``$~\exists ~h~[~\mbox{Subst}    (    g    ,    h    )~\wedge~
\mbox{ExPrf}_{UNION(\beta)}(    h    ,    t    ,    p   )\, ] \, \,$''$~$,
  but \el{encode} is 
 a $\Delta_0^Q$ 
formula --- {\it unlike} the quoted
expression.)
\begin{equation}
\label{encode}
\mbox{Prf}_{~UNION(\beta)}~( \, t \, , \, p \,)~~~\vee~~~\exists ~h\leq p
~~[~ \mbox{Subst} \, ( \, g \, , \, h \, )~\wedge~
\mbox{ExPrf}_{~UNION(\beta)}~( \, h \, , \, t \, , \, p \,)~]  
\end{equation}

% \bbl

Utilizing \eq{encode}'s particular
 $\Delta_0^Q$ 
encoding for 
$\mbox{SubstPrf}_{UNION(\beta)}(   g   ,   t   ,   p  )$, it is
easy to encode
IQFS$(\beta)$'s
 Group-3 axiom, specified by 
\dfx{def-3.12}'s
statement $\oplus \, \oplus \, \oplus  $,
 as  a $\Pi_1^Q$ sentence. $~$  
Thus,  let
$~\Gamma(g)~$ 
denote  \el {encode2}'s formula, and let  $~n~$ denote $~\Gamma(g)$'s
G\"{o}del number.
Then  $~$``$~\Gamma(~ \underx{n}~)~$''$~$ 
is a $\Pi_1^Q$ sentence encoding
IQFS$(\beta)$'s
 Group-3 axiom.
\begin{equation}
\label{encode2}
\forall  ~p~~~
\neg ~ \mbox{SubstPrf}_{~UNION(\beta)} \,
 (  \, g \,   , \, \ulcorner \, 0=1 \, \urcorner  \,  , \,   p \,  )
\end{equation}
In other words,
 $\oplus \, \oplus \, \oplus  $ is encoded as:
``$\forall  ~p~
\neg ~ \mbox{SubstPrf}_{~UNION(\beta)} \,
 (  \, \underx{n} \,   , \, \ulcorner  0=1  \urcorner  \,  , \,   p \,  )$''.

\smallskip

$\Box$

\begin{remm}
\rm
We shall not provide IQFS$(\beta)$'s Groups 1 and 2 axioms with
$\Pi_1^Q$  encodings, similar to Lemma \ref{lem6.2}'s construction. 
% in this chapter.
Such is unnecessary because 
these Group-1 and Group-2 axioms
have obvious encodings, unlike the preceding encoding
of IQFS$(\beta)$'s Group-3  axiom,
described above.
\end{remm}

%\smallskip

\subsection{Main Conjecture and Main  Theorem}

%6666622222

\label{6x.2}

% \smallskip

This section will introduce our main conjecture and
accompanying theorem.
The attached appendix
 will
% also 
provide  evidence,
very strongly suggesting
that 
 \cjx{con6.6}
is correct. This evidence 
shall
almost reach  the critical-mass 
level of
 amounting to a proof. We will begin our discussion
by introducing one preliminary definition,
that will help facilitate the relationship
between  \cjx{con6.6}
and Theorem \ref{thmain}.

% \smallskip

\begin{deff}
\label{def6.5}
\rm
A $\Pi_1^Q$ sentence 
$~\forall \, x ~\phi(x)~$ will be said to be
a {\bf Size K Breaking Pont} iff $~\phi(K)$ is false and 
$~\phi(x)$ is true for all $~ x < K~$.
\end{deff}

\newpage

\begin{coj}
\label{con6.6}
% \label{con6.6}
Suppose $~\gamma~$ 
is an axiom system that includes:
\bed
\rm
\item[  A. ]
All IQFS's Group-zero and Group-1 axioms (from \dfx{def-3.12} ), 
\item[  B. ]
A set of 
 additional
$\Pi^Q_1$ axiom statements, all of which hold
true in the Standard Model, {\it except for one
special
 $\Pi_1^Q$ sentence}, called  $  \Psi  $,
% which
that
 constitutes  a   Size $K$ Breaking Pont
(with $K>2).$
\ennd
Suppose $~P~$ is a proof from  $~\gamma~$ of $~0=1~$ 
(using again \dfx{def-pain}'s $d_{ER}$ deductive method).
Then $~P~$, viewed as a G\"{o}del number, will impose the
following constraint on $~K\,$'s magnitude
\beq
\label{happy}
\mbox{\rm Log}_2 \, K
 ~~~~ <~~~~~ \frac{1}{6} ~~\mbox{\rm Log}_2 \, P 
\enq
\end{coj}
 
We are essentially 100 \% confident that \cjx{con6.6}
is true. Indeed, \el{happy}'s estimate of an upper bound on
$K$'s length is actually an
excessively conservative over-estimate.
An 
approximate
intuitive justification of this inequality is provided in the
attached appendix.
(It actually 
falls 
only
one 
tiny
%small
iota
% away from
short of
 being a 
formal proof.)
Our immediate goal is to discuss and prove
the  \thx{thmain}.
% , which is actually a direct consequence of \cjx{con6.6}.
We recommend that the reader 
postpone examining the Appendix's justification of
\cjx{con6.6}, until after 
 \thx{thmain} and its proof 
%are first
% has been RETURN TO OLD ABOVE 
have been both first read.

%thoroughly examined.

% 
% , before turning to the Appendix's
% explanation of why it is almost 100 /%
% certain that \cjx{con6.6}) is true..)

\begin{theorem}
\label{thmain}
The 
%% formal
\cjx{con6.6}
% implies 
 does imply 
that
$~$IQFS's$~$  axiomatic framework  
satisfies
Definition \ref{def-3.9x10}'s
Consistency Preservation property.
(E.g. that if $\beta$ is an axiom system consistent with
IQFS's Group-0 and 1 axioms then 
% IQFS$(~\beta ~)~$
IQFS$(\beta )$
 is also consistent.) 
\end{theorem}

{\bf Proof.} Similar to the Consistency Preservation 
analysis in our earlier papers, the  proof of 
 \thx{thmain} 
will be a proof by contradiction.

Let
 UNION($\beta$)
once again
  denote the union of 
IQFS$(\beta)$'s Group-Zero,
Group-1 and Group-2 axioms.
If \thx{thmain} was false then it easily follows that
all the axioms in
 UNION($\beta$) will hold true in the standard model,
but  IQFS($\beta$) 
will be contradictory. We will show this is
impossible.

In particular, if  IQFS($\beta$) was contradictory, then there must
exist a minimal sized proof of $0=1$ from IQFS($\beta$), called say $P$.
We call the reader's attention to the fact that 
 UNION($\beta$) consists of all of
IQFS($\beta$)'s axioms,
 except for its Group-3 axiom.
\lem{lem6.3} showed
that IQFS($\beta$)'s Group-3 axiom
can be encoded as
a $\Pi_1^Q$ sentence of the form:
\beq
\label{smile-frown}
\forall ~x ~~~~ \psi(x)
\enq
Then since $P$ denotes the minimal proof of  $~0=1~$ and
since IQFS($\beta$)'s Group-3 axiom states
{\it ``There exists no proof of  $~0=1~$ from
IQFS($\beta$)''}, it follows that
this italicized sentence is false. Thus,
 \el{smile-frown} 
must own a
% owns some 
Size$-K$ Breaking Point, for some $K$.

The \cjx{con6.6}
% explicitly 
indicates
% that
 $K<P$ (because of
\eq{happy}'s inequality). This is impossible because
$P$ was defined to be the minimal proof of  $~0=1~$, and
$K$ represents a smaller
such
 proof of $~0=1~$.
Hence, the
IQFS($~ \bullet     ~$) framework must be consistency-preserving to
avoid
this
% such a 
contradiction.
$~~\Box$

\begin{remm}
\rm
Our earlier  Consistency  Preservation results in 
 \cite{ww93,ww1,ww5,ww9} also relied upon proofs-by-contradiction, mostly
analogous to \thx{thmain}'s verification. 
One
important
%difference,
distinction is,
 however, 
that \thx{thmain}'s proof employed 
\cjx{con6.6}'s paradigm as an intermediate step.
There are good reasons for believing this conjecture is
correct. In particular, the attached appendix provides strong
evidence confirming this conjecture, in a context where its informal
justification falls 
only
one 
small
iota short of being a fully
formalized
 proof. 
At this juncture, we encourage the reader to 
take a quick
glance at
the Appendix's justification for \cjx{con6.6}.
\end{remm}

\section{Further Results and Useful Added Perspectives}

\label{ss7}

Both the  virtues and 
% possible
drawbacks of the IQFS formalism
are consequences of
 Proposition \ref{th-3.3}'s 
characterization of
the  $O\{~[~$Log$(n)~]^3~\}$ 
%invariant for 
quantity of logical symbols 
used for encoding
% needed to to encode 
an integer $~n~$ as a
 grounded term
 $T_n$. 
Thus,
this
% amount 
quantity
is clearly
%much
significantly
 better than the alternate
$O(~n^2~)$ 
length
that arises when the $~\theta~$ function symbol
is replaced by the less efficient primitive 
$~\glamb~,~$ 
%%%%% operator 
defined by Lines 
\eq{zm1} and \eq{zm2}.
On the other hand, one would ideally 
prefer our ground terms to resemble
% the 
conventional encodings
of a binary number that use 
 $O\{~$Log$(n)~\}$ 
logical symbols to encode an
arbitrary  number $~n~$.

It turns out  it is possible to improve 
Proposition \ref{th-3.3}'s encodings 
to such a compressed  $O\{~$Log$(n)~\}$ size,
if one adds only a minor
wiggle to the logic's notation
%%%%%%%%%%%%%%%%%% conventions 
convention.
This distinction arises because most traditional
logic
languages
% , as typically formalized in textbooks,
%% will
formalize
% both conventional terms and 
Definition \ref{def-3.4}'s 
``Ground terms''
as tree-like structures.
An 
easy
alternative
mechanism will 
%%  
%% modification of this 
%% %perspective 
%% construct
%% shall
%% %will
%% % would 
%% 
allow these terms to
own the more
% general 
 generalized
structure of
a Directed Acyclic Graph (Dag).
 
% Kozen has noted ????.
%some computer scientists have noted.

Traditionally, 
this distinction has been
%traditionally 
viewed
 as an unimportant 
wrinkle
%issue 
because it can be 
%readily
easily
proven that 
almost
every Dag-oriented  term can be 
converted into its Tree-oriented counterpart with 
%merely a
  usually only a
% minor
% only an usually unimportant  
Polynomial increase in 
length.
%space. 
The reason for our special interest in
% an alternate 
Dag-formulated ground terms
is that the IQFS formalism 
% has 
possesses
access only to the
specialized
$\, \theta \,$ primitive for generating growth
among integers.
% 
% %distinction
% %between a Tree-oriented and Dag-oriented
% base-language for logic 
% is that the latter 
% ushers in
%  more efficiently encoded Ground terms. 
% 
Thus,   Proposition \ref{th-7.1} 
will
indicate
that
Proposition \ref{th-3.3}'s earlier
 ground terms can 
have their lengths
nicely
 compressed from an
$O\{~[~$Log$(n)~]^3~\}$
 size 
into 
a more 
enticing
%     desirable
% an 
  $O\{~$Log$(n)~\}$ 
magnitude,
when a Dag notation is employed.

%%  in
%% such
%%  a
%% particular
%%  Dag context.

%in the latter context.

\nop
%\newpage

\begin{propp}
\label{th-7.1}
% \rm
Let us consider a Dag-analog of 
\textsection \ref{ss4}'s formalism where once again:
\bee
 \baselineskip = 1.15\normalbaselineskip
\item
$\theta$ is the only available growth permitting function symbol,
\item
the only built-in constant symbols are, again,
% the entities 
 $C_0$, $C_1$ and $C_2$ , 
for representing
%the values of 
0, 1 and 2,
\item
and the three function symbols for
representing
 integer-subtraction, integer-division
and the maximum operation 
are, once 
again, available. 
\ene
%Within this notation, 
In this context,
any integer $~n~$ can be encoded
by a Dag-oriented Ground term 
$~G_n~$
using only
 $O\{~$Log$(n)~\}$ logical 
symbols. { (As the pointers needed to
separate these  $O\{~$Log$(n)~\}$ logical objects 
% will 
use
%have encodings using 
 $O\{~$LogLog$(n)~\}$ bits, the total amount of 
memory 
required
for encoding a
 Dag-oriented Ground term
will employ
at most
% will require
 $O\{~$Log$(n)~\cdot ~~$LogLog$(n)~\}$ bits.)}
\end{propp}

The proof for Proposition \ref{th-7.1}
rests essentially on a more elaborate version of 
\textsection \ref{ss4}'s
justification of  Proposition  \ref{th-3.3}.
It
% essentially rests on 
uses
the fact that 
Proposition  \ref{th-3.3}'s ground term $T_n$
do have
% has
many repeating subterms that can be compressed into
% one 
single objects under a Dag-style notation.

%% 
%% A  $\, 1 ~\frac{1}{2} \,$ page description, of how
%% exactly such a proof of  Proposition \ref{th-7.1}
%% is formulated, is given below:
%% 
%% {\bf Detailed Justification of  Proposition \ref{th-7.1}:} It is easy to

\medskip

{\bf Detailed Proof:} It is easy to
convert the preceding paragraph's summary of the intuition behind
 Proposition \ref{th-7.1} into a formal proof.
Let 
%fffff
$~G_n~$
  denote our Directed Acyclic Graph (Dag)
 for representing an integer $n$,
and let
 $~M~$  be an abbreviation for the quantity
$~\lceil~1\, + \, $Log$_2(n)~\rceil~$.
Our directed graph
$~G_n~$
 will consist of approximately
$~5 \, \cdot \, $Log$(n)~$ nodes.
The first five of its six groups of nodes
in $G_n$'s graph
are defined below 
in 
% xx n roughly
% bottom-to-top order:
bottom-up order:
\bee
\tttc
\item The bottom-most nodes in $G_n$'s
graph
will correspond to the three  
built-in constant symbols of
 $~C_0~$, $~C_1~$ and $~C_2~$, 
that represent
the values of 0, 1 and 2.
%\ene
%\end{document}
\item
Let
 $~\zzthe^j(x)~$
denote the term
 $~\zzthe(~\zzthe(~ ... \zzthe(x)))~$
where there are 
$~j~$ iterations of the 
 $~\zzthe~$ operation.
For each $~j \leq M\,$, the next $j$ levels of
$G_n$'s directed graph will define 
nodes $A_j$ that formalize the quantity
 $~\zzthe^j(1)~$. In a context where 
$~A_0~$ 
is an abbreviation for $~C_1~$,
the remaining $A_j$ are defined by:
\beq
A_j ~~ =~~~ \zzthe(~A_{j-1}~)
\enq
\item
For each $~j \leq M\,$, let $B_j$
denote the value of Max$(A_0,A_1,A_2,...,A_j)$.
In a context where 
$~B_0~$ 
is an abbreviation for the entity $~C_1~$,
the remaining $B_j$ are defined
chronologically in $G_n$'s directed graph by:
\beq
B_j ~~ =~~~\mbox{Max}(A_j,B_{j-1})
\enq
\item
For each $~j \leq M\,$, let $D_j$
denote the value of $~2^{-j} \, \cdot \, B_M~$.
In a context where  
$~D_0~$ 
is equivalent to
%was defined by 
the prior entry $B_M$ in
our directed graph, 
the remaining $D_j$ nodes in our graph
will be defined via \eq{usediv}'s 
Division operation:
\beq
\label{usediv}
D_j ~~ =~~\frac{D_{j-1}}{2} ~~~~~~\mbox{e.g.}~~~~~~ 
D_j ~~ =~~\frac{D_{j-1}}{C_2}
\enq
\item
For each $~j \leq M\,$, 
the node
 $E_j$
will represent the quantity
$~2^j~$.
These nodes in $G_n$'s graph will 
be defined by 
\eq{usediv2}'s 
Division operation:
\beq
\label{usediv2}
E_j ~~ =~~\frac{D_M}{D_{M-j}}
\enq
\ene

Some added notation is needed to describe the last part of
$G_n$'s graph for formalizing $n$'s representation as a
Dag-oriented ground term
employing $O \{ ~$Log$(n)~\}$
logic symbols.
 Let $T_n$
 denote 
Proposition \ref{th-3.3}'s formulation of $~n~$ 
as a Tree-oriented ground term, and 
$G_n$
denote its Dag counterpart
(using the five intermediate steps
itemized above).
Our 
prior
proof of 
Proposition \ref{th-3.3} noted 
$T_n$
% could be 
was
constructed by setting $E_M$ equal to
the least power of 2 greater than $~n~$ and
then  subtracting from it
those powers of 2 which are needed to produce the quantity $n$.

The exact same methodology will now
be used to construct
our 
$G_n$
 representation of $~n~,~$ except 
 we
will now
obviously
use the methodologies from Items 1-5 to
assure
%  that 
no more than $O\{~$Log$(n)~\}$
graph nodes are used to construct
all 
\el{usediv2}'s 
 $E_j$ terms.
(For example since $86~=~128-32-8-2$  
%% $118~=~128-8-2$
which in turn equals ``$~E_7-E_5-E_3-E_1~$'',
 the final stage of 
% $G_n$'s
our
% analogous
 construction of $G_{86}$ will
first
 set node
$F_1$ equal to ``$~E_7-E_5~$'',
then
 set node
$F_2$ equal to 
 ``$~F_1-E_3~$'' 
and lastly have the 
% desired
output node
$F_3$ represent the final answer as the
quantity of  ``$~F_2-E_1~$''.)

\gvx

It is easy to see that this methodology will never use more
than  
  $O\{~$Log$(n)~\}$ logical symbols to encode
 $G_n$ as a Dag-oriented ground term.
Moreover, the needed pointers in the Dag graph $G_n$
will require no more than LogLog$(n)$ bits to
%separate its
distinguish between these
   $O\{~$Log$(n)~\}$ 
separate objects.
%logical symbols.
Hence
if one selects to use a pointer methodology to formulate
$G_n$'s graph,
then
%   our full graph $G_n$ will need 
no more than
 $O\{~$Log$(n)~\cdot~$LogLog$(n)~\}$ bits
will be needed
to encode
all these pointers
(as the last sentence of Proposition \ref{th-7.1} 
had claimed).
  $~~\Box$

% \nyp

\begin{remm}
\label{rem-7.2}. \rm
Let  IQFS$^*$ 
denote the analog of our  
IQFS framework that has   
Proposition \ref{th-7.1}'s
 Dag-oriented Ground Terms replace
  \textsection \ref{ss4}'s earlier 
Tree-oriented Ground Terms
in an accordingly revised language.
A construction exactly analogous to Theorem \ref{thmain}
% 's proof analysis 
will show 
 IQFS$^*$  
also satisfies the
Consistency Preserving paradigm, assuming
(as we do again) that
\cjx{con6.6} is correct.
 Thus,  IQFS$^*$   is a better
formalism, although 
%a bit
somewhat
 more complicated to describe.
\end{remm}

%\medskip

\begin{remm}
\label{rem-7.3}.
\rm
It  
should
 be mentioned that the infinite
number of axiom
sentences,
 appearing in the Group-2 schemas
for ISCE$(\beta)$, IQFS$(\beta)$  and IQFS$^*(\beta)$,
can be
% should be able to be
nicely
 reduced to a 
purely
finite sizes, with almost no loss 
in
%of useful
information. This was done in \cite{ww14} 
for the Group-2
scheme
of its IS$_D(\beta)$ formalism, 
with the latter
%%%%%%%%%%%%%%%% still 
%where the
% 
% germane Group-2 scheme was
% reduced to one
% single  axiom sentence 
% while the resulting 
% 
% latter
%formalism still
% produced
producing
isomorphic counterparts
of all of $~\beta \,$'s
full set of
 $\Pi_1$ theorems
(e.g. see Sections 5 and 6 of    \cite{ww14}).
The same methods will
% trivially 
routinely
generalize for the
% 
% Analogs of the techniques from Sections 5 and 6 of 
% \cite{ww14} 
% % will easily 
% apply to each of the 
% 
ISCE,
 IQFS and IQFS$^*$ frameworks.
It will imply that easy modifications of
IQFS(PA+) and IQFS$^*$(PA+), 
owning a 
{\it strictly finite number} of axiom sentences,
are able to prove isomorphic counterparts of
{\it the full infinity} of $\Pi_1$ theorems, produced by Peano Arithmetic 
\end{remm}

\begin{remm}
\label{rem-7.4}
\rm
There is one other amendment to the
 IQFS and IQFS$^*$ 
formalisms that should be mentioned because of
its pragmatic significance.
% from an engineering perspective. 
Let $S_j$ denote
\eq{eng}'s sentence. Let  
 IQFS$_R$ and IQFS$_R^*$ denote the minor modifications
of   IQFS and IQFS$^*$  that include $~S_1,~~S_2,~~S_3,~...~$ as
%additional 
axioms. 
\beq
\label{eng}
E_j~-~E_{j-1}~~~ = ~~~E_{j-1}~~~
\enq
The addition of these $S_j$ as axioms
 cannot possibly affect IQFS's  consistency
because each  $S_j$ is provable
from the Group-Zero and  Group-1
% and Group-2 
axioms as a theorem.
The advantage of treating such  $S_j$ as
such built-in axioms
(rather than as theorems)
is that many of IQFS$(\beta)$'s other
theorems can have their proofs conveniently shortened in 
length
in this case. 
\end{remm}
%% 
%% The advantage of treating the  $S_j$ as axioms
%% is that many of IQFS$(\beta)$'s other
%% theorems can have their proofs  undergo an 
%% exponential shortening of their proof lengths, in this case. 

\begin{remm}
\label{rem-7.5}
\rm
Our
 proposed
$\theta$ primitive is
likely to be
  of philosophical interest,
quite apart from its implications
for
Hilbert's
Second Open Question. This is because:
% primary
%main
%  two
%  Some reasons for this are summarized below:
\bee
%\small
\baselineskip = 1.24\normalbaselineskip
%  \baselineskip = 1.2\normalbaselineskip
\parskip 0 pt
\item It is interesting that
%  
% Proposition  \ref{th-7.1} showed 
% our newly proposed
% % that the
% %%  formally
% ``indeterminate  and 
% growth-permitting''
% 
the $\theta$
primitive, combined with 
the 
% traditional
%contrasting
{\it non-growth}
subtraction, division and maximum operands, can
represent any integer $n$  with
% essentially 
a complexity comparable to 
the $O(~$Log$(n)~)$ 
sized
representations of  an integer 
under its
 conventional      binary encoding.
{\it $~$Does this
% mean 
suggest the 
%primitive 
foundational
reasoning capacities of either a child
or 
of the anthropological predecessors of modern man
% had 
%did 
possess
access to
a technique for reasoning that
% could 
constructs the natural numbers
without the traditional
uses of the
 addition and  multiplication
functions?}
%primitives?} 
There is likely to  be no easy
 yes-or-no answer
 to this
question.
(Our suspicion is that an
infant-child's
 primitive thoughts
rely upon an
analog of 
 $\theta$'s basic indeterminate
functionality.)
% function.)
\item
Also, it is philosophically 
% tantalizing 
 curious
%interesting
whether or not
the Group-2 axioms of 
Definition \ref{def-3.12}'s
IQFS formalism
% manage to
capture the
majority of the engineering implications of Peano Arithmetic?
This is
because 
one may
% thus
% consequently
 wonder whether
the {\it purist engineering implications} of Peano Arithmetic
%do 
consist of
its particular
% especially 
% special  
%derived 
theorems,
% that can be formally
 that can be
 encoded
%% , very formally,
%under
in such
% a 
%%%%%%  precise
%particular
%special
 ``$~\Pi_1^Q~$formats''$~?$
%(perhaps conjoined with Remark \ref{rem-7.3}'s compression
%methodologies)$~$?
\end{enumerate}
There is no space available in this 
% short
article to
% answer
adequately address 
%%or
%even 
%%%merely
%delve, with 
%in adequate detail,  
these
% preceding
philosophical
issues.
%questions.
Our point is, however, that IQFS
raises other issues 
besides
%our 
besides this article's
%central 
 main
theme,
%  the central issue 
that
% {\it some 
% carefully
% selected fragments}
some  fragments
 of the 
historic
% initial
% particular
aspirations of
Hilbert and G\"{o}del
 in 
%raised in
  $*$ and $**$ should be revisited.
\end{remm}
% \gv2

\section{Disentangling Some 
%Incorrect 
Confusing
Interpretations
 of 
Theorem \ref{thmain} and Its Corollaries
that
 Could
% Might
 Otherwise
Plausibly
% Easily
%Potentially 
 Arise}

\label{ss88}

% hhh8888
It was
during one of the final stages of this paper,
as its draft was nearing 
completion, that we elected 
to add a very short
further
% $1~\frac{1}{2}$ page
chapter into this article, 
{\it ``disentangling''}
some
% serious 
% levels of
 confusions that might 
{\it ``otherwise
% `could
plausibly arise''} .
It was
desirable
 to 
insert a
% an added 
% brief
chapter into this article,
 whose title
contained
the preceeding italicized words,
so that the 
{\it exact meaning}
 of 
% its 
Theorem \ref{thmain}'s central result
 could not
% cannot
%% be 
% possibly
 be 
 confused.

%% 
%% Theorem \ref{thmain}'s
%% % exact
%% % Its  
%% implications 
%%  are 
%% %also 
%% subtle
%% %partially
%% partly
%% because it is
%% somewhat
%% % quite
%% %obviously
%% unconventional for an author$~ X~$ to reply to an
%% % non-trivial
%% open question, raised by a
%%  world-famous mathematician $~H~,~$ by 
%%  first editing
%% % the nature of
%%  the {\it very question}
%% XYZ that $~H~$ had
%% initially
%%  raised.
%% Yet,
%% % precisely 
%% roughly
%% this  was  done in  
%% % what we had done in  
%% \textsection \ref{ss4}, when we
%% had divided Hilbert's Second
%% Open Question into two
%% quite
%%  separate problems, which were called the
%% Q-1 and Q-2  variations of Hilbert's historic
%% % millenial
%% Second
%%  Open Problem.
%% 
%The meaning of  Theorem \ref{thmain} is
%also

Theorem \ref{thmain}'s meaning is
% further
 complicated by the
fact\footnote{$~~$PA+ was 
defined in
% several of
our  articles to be the trivial extension
of Peano Arithmetic that uses a language that includes the
six Grounding level functions, in addition to the usual
addition and multiplication primitives.}
 that 
%our 
IQFS$(PA+)$,
% formalism, 
as well as its 
several
%formalized 
predecessors
in
%the article 
\cite{ww93}-\cite{ww14},
were able to prove 
% a
theorems,
%formally
 asserting
 their own 
consistency,
only by employing
Example \ref{ex-2.5}'s 
built-in
{\it ``I am consistent''} statements,
as 
% a formalized
invoked 
 axiomatic 
declarations.
% statement Such an axiom, 
This approach,
 which 
was
%had been
centered around 
%the statement 
Item                       $~\oplus~$
from Example \ref{ex-2.5},
will naturally cause many
readers to
%  ask
raise
 the following
%  skeptical-but-quite-legitimate
%partially  skeptical 
%and quite reasonable
skeptical
%subsequent
question:
\begin{quote}
$\# \, \# \, \#~ $
Is it not {\it almost cheating} when an axiom system verifies
its own consistency by using
%  statement 
 $~\oplus \,$'s formalized
{\it ``I am consistent''} axiom as an intermediate step, to verify its
own consistency? After all, such a
% formalism 
technique
can verify its
own consistency
only in
a {\it technically purely
legalistic sense}.
(It is
%  albeit
certainly,
however,
 not
meaningful
 in the much {\it broader almost
philosophical
% sense,}
respect,}
 that Hilbert was grasping for in his
famous 
year-1900
Second Open Question.)
\end{quote}
It will be necessary
 for us 
to introduce
the formalisms of
Definition \ref{def8.1} 
and Corollary \ref{cor8.2},
before replying
 to
~$\# \, \#  \, \# \, $'s 
%important 
%open
%particular
%exact
%%% precise
%important
 query.
%  question.

\begin{dff}
\label{def8.1} 
\rm
$~$ 
Let $~\beta~$
% once 
again denote an axiom system that
uses $~L^Q~$'s language (whose
collection of 
 function symbols include the
six
 usual
% six
 Grounding-level primitives plus $\theta$'s special
indeterminately-defined function operation).
Also, let IQFS$(\beta)$ 
once
again denote
Definition \ref{def-3.12}'s
% \textsection \ref{ss5}'s proposed 
self-referencing  formalism, that proves
all of $~\beta\,$'s  
 $\Pi_1^Q$ theorems
 and which 
additionally 
% can
recognizes its 
own formalized
 Hilbert
%-style
consistency,
via
its use of an
invoked
% Definition \ref{def-3.12}'s
% Group-3 
{\it ``I am consistent''}
axiom. 
%% 
%% % (through 
%% (via 
%% again
%% the use of 
%%  Example \ref{ex-2.5}'s
%%  {\it ``I am consistent''}
%% % axiomatic
%% axiom).
%% %  statement). 
%% % the formalism
 Then  $~\beta~$ will be
%defined to be
called
  {\bf Platonically Stable} iff   
   IQFS$(\beta)$ is a
% formally 
consistent axiom system.
% Likewise,
( Likewise, if IQFS$(\beta)$ is inconsistent then  $~\beta~$
will be called {\it ``Platonically Unstable''}.)
\end{dff}

%% The arithmetical system
%%   $~\beta~$
%%  will
%% % shall
%%  be
%% called {\bf Platonically Unstable}, otherwise.
%% 
%% 
%% {\bf Comment:} If IQFS$(\beta)$ is inconsistent then likewise $~\beta~$
%% will be called {\bf ``Platonically Unstable''}.

\begin{ccr}
\label{cor8.2}
\rm
$~$ 
Let us
% again 
assume
% that  
\cjx{con6.6} is correct
(as we are almost 100 \% certain
% that 
it is).
%Also,
Suppose 
 $\beta$ is 
an axiom system which is consistent with
IQFS's Group-0 and 1 axioms.
%and that \cjx{con6.6} is correct.
 Then
 $\beta$ is automatically {\it ``Platonically Stable''}
(e.g. the formalism IQFS$(\beta )~$ is
% thereby,
% thus,
%automatically 
consistent). 
\end{ccr}

The formal statements of 
Theorem \ref{thmain} and
 Corollary \ref{cor8.2}
 are
closely related. 
%nearly identical.
This is because
 Corollary \ref{cor8.2} applies
Definition \ref{def8.1}'s notion of  ``Platonic Stability''
in essentially
 every
% each
 place where Theorem \ref{thmain} had used
Definition \ref{def-3.9x10}'s ``Consistency Preserving'' paradigm.
Thus,
there is
% , thus,
  no need to
 formally  prove  Corollary \ref{cor8.2}
%(since it is,
because it is
% obviously,
% a trivial 
almost a direct 
%an immediate
consequence of  
Theorem \ref{thmain}.

% 's underlying  paradigm.

\medskip

The 
nice aspect of
 Corollary \ref{cor8.2}
is
% that
 its  vocabulary
 will 
%shall
 allow us to 
better
% more easily 
frame a
% an exact response
%  much  clear  
reply  to 
$~\# \, \#  \, \# \, $'s 
% pressing
% important
query.
% query. 
In particular, it would
% certainly
% be 
 certainly 
be preferable (and more
% ideally
Utopian) if human beings could muster
% a 
more 
sophisticated
% complicated 
justifications
for their thought processes than had appeared in  
 Example \ref{ex-2.5}'s
short
 1-sentence axiom statement $~\oplus~$
(which
% had abruptly 
 simply  declared
%  merely
 {\it ``I am consistent''}).
However, 
such
a
% more elaborate
%detailed
reply to $~\# \, \#  \, \# \, $'s 
% legitimate 
%reasonable
query
is
% actually
%strictly
  unnecessary for  the
planet Earth's
% highest-IQ 
 high-IQ 
primates
% specie, residing on the planet Earth,
 to gain
% an 
adequate 
%threshold
confidence
 to
justify their
% introspective
 thought processes. All that is 
 {\it strictly} necessary, {\it  from a 
% technically
purist perspective,} is for an advanced
Thinking Being   
to find a
robust formalism that satisfies Definition \ref{def8.1}'s 
criteria of  ``Platonic Stability''. 
The 
%nice 
significant
aspect of our
 Corollary \ref{cor8.2} is 
 it shows
%ssssssss
adequate 
forms of
 Platonic Stability
are available to
%show
allow 
an introspective thinker to simultaneously:
\bee
\small
\item
presume its own consistency as a built-in assumption, and
\item
rest assured
%% that 
this  assumption will not spin
its
%system  
IQFS$(\beta )$
formalism 
into a cycle of inconsistency.
\ene
%One nice aspect of
%  Corollary \ref{cor8.2} hits
Corollary \ref{cor8.2}'s notion of  Platonic Stability, thus,
% may
partially reinforces
G\"{o}del's philosophy of
Mathematical Platonism,
and it also
%%  seems to
helps
make comprehensible
%  {\it at least  fragments} of 
the aspirations
that motivated Hilbert's and G\"{o}del's 
famous statements $*$ and $**$.
% , that had appealed to G\"{o}del and many other scholars.
%Moreover, 
This
%its
 perspective is 
% also 
useful
%%% 
%%% % Moreover, it is useful 
%%% Moreover, its
%%%  perspective is useful
when we remind ourselves
% remind
% our readers
 that the invariant $~++~$, which 
Example \ref{ex-2.3}
had attributed to the joint work
of  Pudl\'{a}k, Solovay, Nelson and Wilkie-Paris,
\cite{Ne86,Pu85,So94,WP87},
%will show 
established
that self-justifying formalisms, much stronger than
  IQFS$(\beta )$, are
% also
  {\it simply impossible.}

% hh9999

\begin{remm}
\label{rem8.3}
% {\bf (repeating our reply to  $~\# \, \#  \, \# \, $'s query} 
%{\bf (briefly repeating our main theme):}
{\bf (Snapshot Perspective):}
\rm
Summarizing  our last
25 years of research
into one short paragraph,
it is
certainly
 true that any  proof that relies
upon 
 Example \ref{ex-2.5}'s
{\it ``I$~$am$~$consistent''} 
axiom
is, in some respects, 
a quite
%% {\it very}
$\,$skinny
%slender
form
%% variant
of proof, that
one is 
{\it   almost first tempted}
to$~$ignore. On the other hand,
this perspective  is  actually
quite useful
%  helpful
because    Corollary \ref{cor8.2} 
does reminds us
that
we 
do live
%reside
in 
what Definition \ref{def8.1} 
%  had 
had called a
{\it  `` Platonically Stable''}
world. 
\end{remm}

%%%%%bbbbbotom
%%%%%bbbbbotom

\section{Concluding Remarks}

%    
%\vspace*{- 0.4 em}
%vvvvvv

All  our published articles
% (since 1993) 
 about self-justifying arithmetics
have emphasized 
that
% our proposed
 evasions of the
Second Incompleteness Effect  rested
% upon 
on
using arithmetics
that were
% much 
weaker than traditional arithmetics in, at least,
some 
% particular
%well-defined 
respects. 
%(Thus, the  
(The
% importance 
% % significance
% of the 
Second Incompleteness Theorem's
%important 
%underlying
% vital 
significance
in refuting 
%ss% % implications
%ss% in establishing a
%ss% 90-95 \%
%ss% refutation of 
the original objectives of
Hilbert's Consistency Program
% are, 
is
thus,
simply,
% certainly, 
%of course, 
%%%%% thus,
%undeniable. 
%certainly 
% unquestionably 
undeniable.)
%omni-present.
It would, nevertheless, be of interest if
some
%% 5-10 \%
fragments of 
%its
Hilbert's
and G\"{o}del's 
 objectives, as
%% was
 announced in
their  formal
 statements
$*$ and $**$, 
were
% actually
%partially
 achieved. 

\smallskip

It is in this context and 
also
 with knowledge about how the scope
of the Second Incompleteness Theorem has been
%
% significantly
% greatly
% expanded
enhanced
 in
%several
many
 respects
% in 
during
the last 35 years, that we hope 
% that
our boundary-case exceptions to the Second Incompleteness paradigm
%will
shall
 be found 
to be
useful.

%% may 
%% % still 
%% prove 
%% to be 
%% useful.
%% %effective. 
%interesting.

%useful.

%\smallskip

This is not merely because it would be
pleasantly
 reassuring
to see
achieved
 some 
select
%particular
portions of what early 20-th century articles
had advocated.
(For instance, G\"{o}del's 
%documented 
preference \cite{Da97,Go5} 
for a ``Platonist'' philosophy 
%% 
%% for mathematics 
%% has been
%% %nicely 
%% sharply
% was partially
is partially
%has been 
reinforced by 
 Corollary \ref{cor8.2}'s 
% demonstrated
%  of 
% some useful
%  a
 ``Platonic Stability''.)
Our
% research is, 
results are,
also  of interest
% significant
% interest 
%  It is also 
because 
future
 computers 
will,
% of the future may,
 perhaps,
% will hopefully
have their capacities 
%%%significantly
%advanced, 
% increased,
enhanced,
if they can simulate
a human's
% being's 
gut instinctive faith in his
own consistency.  

\smallskip

Thus, it is in 
such respects
% this respect
that 
Theorem \ref{thmain} and
% also 
the related
Propositions \ref{th-3.3} and  \ref{th-7.1} 
%  and  \ref{th-7.4}
% shall 
will
  help us
%  to
% much 
better
appreciate the aspirations that Hilbert and G\"{o}del
were pursuing in
their formal
 statements
$*$ and $**$,
$~$as well as
 to  formalize
what
%specific
%particular
%fragments
{\it limited}
aspects
of 
their
stated
 goals in $*$ and $**$ can
%  be 
actually
be 
achieved. 

% \smallskip

It
 also 
should   be mentioned
that
 some
special
% interesting
philosophical 
 issues
 were
raised
%both 
by
Remark \ref{rem-7.5}
and
% also 
in
Remark \ref{rem8.3}'s
observation
that symbolic logic
has a {\it  ``Platonically Stable''} structure.

\bigskip
\smallskip

% \medskip

{\bf ACKNOWLEDGMENTS:}
%%%%  several Sections 1-4, 
%\textsection \ref{ss2}, 
I am
% much
%very 
grateful to
%was influenced by an emailed letter from 
Pavel Pudl\'{a}k
for 
sending me an email letter \cite {Pupriv} that  
suggested
I investigate how to apply
% an analog of 
Ajtai's study
\cite{Aj94} of Pigeon-Hole effects 
to
refine my prior results about self-justifying logics.
(The combination of
 Pudl\'{a}k's
%insightful             
suggestion
% \cite {Pupriv}
and our
subsequent
% further 
        distinguishing
between the
$~\glamb   ~$ and $~\theta~$ operators
has led to the 
stipulated
 improvements upon
\cite{wwapal}'s 
%% 
%% \nyp
%% 
%% \noindent
ISCE formalism.)
% I am very grateful to Pudl\'{a}k for making this
% suggestion. 
I also thank Bradley Armour-Garb
and Seth Chaiken
 for
%  many
several helpful
comments about how to
% significantly
 improve 
% the 
this paper's
presentation.
As was mentioned earlier,
%  in the text, 
Sam Buss
should also be credited for partially
%influenced
 influencing
 my research
when,
during a casual lunch at a 1977 meeting of the
Kurt G\"{o}del
Society,
 he suggested
I might look at the 
implications the
Pigeon Hole Principle
had for my 
on-going
research.

% \bigskip

%% 
%% It
%%  also 
%% should   be mentioned
%% that
%% some
%% special
%% % interesting
%% philosophical 
%%  issues
%%  were
%% raised
%% %both 
%% in
%% Remark \ref{rem-7.5}
%% and also in
%% Remark \ref{rem8.3}'s
%% observation
%% that we {\it   actually  do} (!) 
%% live   
%% in what can be called a 
%% {\it  ``Platonically Stable''}
%% world.

%   \newpage

\section*{Appendix Justifying  \cjx{con6.6}} 

This appendix
should be read only after Section \ref{nn6} has been
finished.
It
will have two purposes. They will be to:
\bee
\item
 Outline an approximate intuitive justification for  
 \cjx{con6.6}'s validity.
\item
Explore this conjecture's broader significance.
\ene

% Explain why it is desirable to supplement this hand-waving argument
% with a
% more
%  formalized proof confirming \cjx{con6.6}.

\subsection*{A.1 $~$ An 
Approximate
Intuitive
Sketch 
 Justifying \cjx{con6.6} }

% aaaaaaaa

The underlying intuition supporting
\cjx{con6.6} 
is
% quite 
easy to summarize.
% It is due to 
It stems from
the fact that
 $\theta$ 
is the only growth-permitting function
(or primitive)
in $L^Q$ 's language.
Moreover,
 the definitions of  $\Pi_1^Q$ and 
$\Sigma_1^Q$ sentences
(in Item \ref{def-3.10})
 {\it $~$forbid$~$} 
 $\theta$'s
 appearance 
{\it anywhere} in these
  sentences, 
{\it outside} of the particular places
%{\it outside the case}
 where
 Definition \ref{def-3.3}'s $E_j$
 terms
   use
the  $\theta$ primitive
as an
intermediate device for 
defining
the quantity $2^j$.

These   facts
% seem to  imply 
 seem to  obviously imply that 
a proof can construct
integers  $n>2$ only by applying the  
Up-Walking axioms,
from Lines  \eq{walk1}--\eq{walk4},
%%%%%% in essentially the very direct
%lengthy 
%canonically
in the directly
 canonical and
cumbersome manner.
%In particular, 
That is,
 a proof $\, P \,$
appears to be able to
% can, seemingly, 
%construct  
verify the existence of
an
integer   $n \geq 2^d$ only after it has isolated $d+1$
 distinct terms,
%% which it has determined represent
which correspond to 
$d+1$ {\it distinctly different} 
powers of 2.

%% (analogous to 
%% Line \eq{ej-def}'s many
%% iterated
%%  terms of the form
%%  ``$~\theta^n(1)~$'' ).
 
Let $M_d$ denote a finite-sized  model 
for the
subset 
of  Natural Numbers,
all of which
satisfy the 
 exact
 inequality of
 $~ x \, < \,2^d~$. Also, let $~S~$ denote 
the set
% of 
of $\Pi_1^Q$ axiom sentences that were defined by
 Part (B) of 
\cjx{con6.6}.
If $~d \, = \, \lfloor \,$Log$_2 \,K \, \rfloor$ then 
Part (B) 
%%% states that
implies that
all of the axioms in $~S~$ 
hold true within the model $M_d$.
However, its 
particular
``Breaking Point'' sentence $~\Psi~$
% becomes 
does  become
%% uniquely
false in the 
model $M_{d+1}$, while all the other axioms of $~S~$
continue to
 hold true
within each of the models of
 $M_{d+1}\, , \,M_{d+2}\, , \,M_{d+3}\, , \,... ~$.

%~  M \inf  ~$.

 At an intuitive (somewhat informal) level,
these facts
% strongly
% hint
suggest
 that
the proof $~P~$,
alluded to in \el{happy},
% is able to 
can
establish ``$0=1$''
only after it 
has
first constructed more than 
$ 1 \, +  \, \lfloor \,$Log$_2 \,K \, \rfloor$ distinct powers of
2. (This is because
the sentence $~\Psi~$
ultimately forces a contradiction, in such a context,
because of its
``Size-$K$ Breaking Point'' property.).

%%  \cjx{con6.6}'s proof $~P~$ can verify the {\it contradictory assertion}
%% ``$~0=1~$'' only after it has verified
%% (at least implicitly) the existence of
%% an integer $~K~$ that violates
%% the requirements of $\Psi$'s  
%%  ``Breaking Point'' sentence  of $~\Psi~$.

 Under our rules for encoding  $\, P \,$
(which requires six 
bits to encode each of $P$'s
23
%employed
 logical
symbols),
such a 
lengthy
sequence  of
 $~d+1~$
 distinctly  different
terms will
% certainly 
% easily
cause $P$ to be so large 
that Log$_2(P)$ 
exceeds by at least a factor of 6 the size of
 Log$_2(K)$. 
Indeed, the quantity $\frac{1}{6}$ in 
\el{happy}'s
 upper bound
% equation 
is
actually a conservative overestimate of the
actual bounding quantity. But we will not delve further into
such a  discussion, here, because this section 
was intended to
provide only an
intuitive sketch justifying
 \cjx{con6.6}.

A 1-sentence summary of the intuition behind \cjx{con6.6}
is that its axiom system $\gamma$ is so weak that it cannot
construct a
required integer, larger than $K$, without engaging in a
proof that
directly
 requires more than $6 \, \cdot \, $Log$_2 K \,$ bits.  

\subsection*{A.2 $~$More Details about \cjx{con6.6} }

Our IQFS($\beta$) axiom system 
{\it does not need}
an actual 
 formal proof of
\cjx{con6.6} 
for it to become mathematically 
very
significant.
All it needs, {\it from a
strictly minimalistic perspective,} is for 
the \cjx{con6.6} to 
hold
 true 
under the Standard Model
(even if its
formal
 proof is
% actually
 independent from
the 
exact
axioms of say Peano Arithmetic and/or ZF Set Theory).

A 
% formalized 
rigorous proof of \cjx{con6.6}'s is, obviously,
still desirable. This proof is likely to be
tedious
and 
%%very 
lengthy because
the analog of \cjx{con6.6}'s inequality  
\eq{happy} is,
simply,
 not  true when IQFS's Group-zero axiom
 replaces
our new
 $\theta$ function operation  with
the Successor primitive.
In particular, the Pudl\'{a}k-Solovay
Assertion $++$, from 
Example \ref{ex-2.3}-a, was 
able to generalize the Second Incompleteness
Theorem, 
only
by using the fact that 
analogs of \cjx{con6.6} are false
when Successor replaces $\theta$.

The intuitive reason 
that $\theta$ operates 
so
differently from Successor
% for this divergence 
is that
these two
operations differ in the following three 
significant
% important
respects:
\bee
\item Integers grow
% directly 
%in
at
 a monotonically uniform 
rate
%manner
under the successor operation, while any finite
 sequence 
$\theta(2),~\theta(\, \theta(2))~,~\theta(\, \theta(\, \theta(2))) ~...
  ~\theta^k(2)~$  is allowed to be monotonically decreasing
under the $~\theta~$ primitive.
\item 
% Also, 
The $\theta$ primitive satisfies Definition \ref{def-2.6}'s
property of being a ``Q-function'' with $\aleph_1$ different
possible
% permissible
solution vectors, while there is only one unique  representation of
the successor function.
\item 
The $\Pi_1$ encodings for
axioms about
 the
 3-way arithmetic predicates $Add(x,y,z)$
and  $Mult(x,y,z)$,
such as the associative rules for addition
and multiplication,
 help one characterize the implications
of the Successor function.
But these rules
are 
fully 
irrelevant
% in regards 
to  $\theta$'s behavior.
\ene
Thus,
Items 1-3 make us essentially 100 \%
confident that the $\theta$ and Successor primitives
shall differ
 sufficiently
% different
 for 
 \cjx{con6.6} to apply to the former
primitive, although certainly not
also to its 
% latter
counterpart
under the Successor operation.
%operator.

\end{document}